 \documentclass[12pt,leqno]{amsart}
\newtheorem{dummy}{}[section]
\newtheorem{definition}[dummy]{Definition}
\newtheorem{theorem}[dummy]{Theorem}

\newtheorem{proposition}[dummy]{Proposition}
\newtheorem{lemma}[dummy]{Lemma}
\newtheorem{corollary}[dummy]{Corollary}
\newtheorem{remark}[dummy]{Remark}
\newtheorem{example}[dummy]{Example}

\begin{document}
\bibliographystyle{plain}
\title{Derived Langlands III: PSH algebras and their numerical invariants } 
\author{Victor P. Snaith}
\date{5 May 2020}
\maketitle
 \tableofcontents 
 
 \[  \begin{array}{c}
{\bf  In \ memory \ of \ my \  friend \ John \ Conway \  (1937-2020)}  \\
{\bf  and \ his \ spherical \  blackboard}
 \end{array} \]

 \begin{abstract}
  In \cite{Sn18} and \cite{Sn20} I studied various topics related to an embedding of the category of admissible representations of a locally $p$-adic Lie group $G$  into the derived category of an additive category  called the monomial category of $G$. By \cite{Sn20} this material had crystallised 
sufficiently that it was all about the hyperHecke algebra and its bar resolution. This is not too surprising because the connection between admissible representations and modules over the classical Hecke algebra is a fundamental connection often used by the Langlands Programme professionals (\cite{PD84}; \cite{AJS79}). The hyperHecke algebra generalises the classical Hecke algebra
and in the case of finite groups it was implicit in the stable homotopy results related to the Segal conjecture \cite{MSZ89} - an area of algebraic topology which led to the first explicit formula for Brauer's induction theorem \cite{Sn88}. However, even in the finite group case the hyperHecke algebra never quite got invented, as far as I know, in the ``Segal conjecture years''.

In another, but related, direction Hopf algebras (notably PSH algebras) entered into admissible representation theory of general linear groups in \cite{BZ76}, \cite{BZ77} and \cite{AVZ81}. 

This sequel to \cite{Sn18} and \cite{Sn20} studies some PSH algebras and their numerical invariants, which generalise the epsilon factors of the local Langlands Programme. It also describes a conjectural Hopf algebra structure on the sum of the hyperHecke algebras of products of the general linear groups over a $p$-adic local field or a finite field.
 \end{abstract}

 \section{Introduction}
 
 This essay is a sequel to \cite{Sn18} and \cite{Sn20}. It is a first attempt at clarifying the role of Hopf algebras in their approach towards the representation theory of locally profinite Lie groups.
 
 I am very grateful to Francesco Mezzadri for showing me his (then) conjectural formula for the analogue for symmetric groups of the Kondo-Gauss sums for $GL_{n}{\mathbb F}_{q}$. The formula appears in (\cite{Sn18} Appendix III, \S4) together with my proof that it respects the PSH algebra product. That formula started me
 thinking about numerical invariants for other Hopf (or PSH-) algebras and thence to wondering whether the hyperHecke algebras of general linear groups formed an analogous Hopf algebra.

\section{PSH-algebras over the integers}
\begin{dummy}
\label{1.1}
\begin{em}

A PSH-algebra is a connected, positive self-adjoint Hopf algebra over ${\mathbb Z}$. The notion was introduced in \cite{AVZ81}. Let $R = \oplus_{n \geq 0} \ R_{n}$ be an augmented graded ring over ${\mathbb Z}$ with multiplication 
\[ m : R \otimes R \longrightarrow R. \]
Suppose also that $R$ is connected, which means that there is an augmentation ring homomorphism of the form
\[   \epsilon :  {\mathbb Z} \stackrel{\cong}{\longrightarrow} R_{0} \subset  R .\]

These maps satisfy associativity and unit conditions.

\underline{Associativity:}

\[  m(m \otimes 1) = m(1 \otimes m) : R \otimes R \otimes R \longrightarrow  R .  \]

\underline{Unit:}

\[    m(1 \otimes \epsilon ) = 1 = m(\epsilon \otimes 1);  R  \otimes {\mathbb Z}   \cong R  \cong {\mathbb Z} \otimes R \longrightarrow  R \otimes R \ \longrightarrow R . \]

$R$ is a Hopf algebra if, in addition, there exist comultiplication and counit homomorphisms
\[ m^{*} : R \longrightarrow R \otimes R \]
and
\[ \epsilon^{*} : R  \longrightarrow  {\mathbb Z} \] 
such that

\underline{Hopf} 

$m^{*}$ is a ring homomorphism with respect to the product $(x \otimes y)(x' \otimes y') = xx' \otimes yy'$ on $R \otimes R$ and $\epsilon^{*}$ is a ring homomorphism restricting to an isomorphism on $R_{0}$.
The homomorphism $m$ is a coalgebra homomorphism with respect to $m^{*}$. 

The $m^{*}$ and $\epsilon^{*}$ also satisfy

\underline{Coassociativity:}

\[  (m^{*} \otimes 1)m^{*} = (1 \otimes m^{*})m^{*} : R \longrightarrow R \otimes R \otimes R \longrightarrow  R \otimes R \otimes R \]

\underline{Counit:}

\[    m(1 \otimes \epsilon ) = 1 = m(\epsilon \otimes 1);  R  \otimes {\mathbb Z}   \cong R  \cong {\mathbb Z} \otimes R \longrightarrow  R \otimes R \ \longrightarrow R . \]

$R$ is a cocomutative if

\underline{Cocommutative:}

\[   m^{*} = T \cdot m^{*} : R \longrightarrow R \otimes R \]
where $T(x \otimes y) = y \otimes x$ on $R \otimes R$.

Suppose now that each $R_{n}$ (and hence $R$ by direct-sum of bases) is a free abelian group with a distinguished ${\mathbb Z}$-basis denoted by $\Omega(R_{n})$. Hence $\Omega(R)$ is the disjoint union of the $\Omega(R_{n})$'s. With respect to the choice of basis the positive elements $R^{+}$ of $R$ are defined by
\[ R^{+} = \{ r \in R \ | \  r = \sum \ m_{\omega} \omega , \  m_{\omega} \geq 0, \omega \in \Omega(R) \}.\]
Motivated by the representation theoretic examples the elements of $\Omega(R)$ are called the irreducible elements of $R$ and if $ r = \sum \ m_{\omega} \omega \in R^{+}$ the elements $\omega \in \Omega(R)$ with $m_{\omega} > 0$ are call the irreducible constituents of $r$.

Using the tensor products of basis elements as a basis for $R \otimes R$ we can similarly define $(R \otimes R)^{+}$ and  irreducible constituents etc.

\underline{Positivity:}

$R$ is a positive Hopf algebra if 

\[ m((R \otimes R)^{+}) \subset R^{+}, m^{*}(R^{+}) \subset (R \otimes R)^{+} , \epsilon({\mathbb Z}^{+}) \subset R^{+}, \epsilon^{*}(R^{+}) \subset  {\mathbb Z}^{+}. \]

Define inner products $\langle - , - \rangle$ on $R$, $R \otimes R$ and ${\mathbb Z}$ by requiring the chosen basis ($\Omega({\mathbb Z}) = \{ 1 \}$) to be an orthonormal basis.

A positive Hopf ${\mathbb Z}$-algebra is self-adjoint if

\underline{Self-adjoint:}

$m$ and $m^{*}$ are adjoint to each other and so are $\epsilon$ and $\epsilon^{*}$. That is
\[ \langle m(x \otimes y), z \rangle =  \langle x \otimes y, m^{*}z \rangle  \]
and similarly for $\epsilon, \epsilon^{*}$.

The subgroup of primitive elements $P \subset R$ is given by
\[   P = \{ r \in R \ | \ m^{*}(r) = r \otimes 1 + 1 \otimes r \}  \]
\end{em}
\end{dummy}

\section{The Decomposition Theorem}

Let $\{ R_{\alpha} \ | \  \alpha \in {\mathcal A}  \}$ be a family of PSH algebras. Define the tensor product PSH algebra 
\[  R = \otimes_{ \alpha \in {\mathcal A}}  \  R_{\alpha}  \]
to be the inductive limit of the finite tensor products $ \otimes_{ \alpha \in S}  \  R_{\alpha}$ with $S \subset {\mathcal A}$ a finite subset. Define $\Omega(R)$ to be the disjoint union over finite subsets $S$ of $\prod_{\alpha \in S} \ \Omega(R_{\alpha})$.

The following result of the PSH analogue of a structure theorem for Hopf algebras over the rationals due to Milnor-Moore \cite{MM65}

\begin{theorem}{$_{}$ }
\label{2.1}
\begin{em}

Any PSH algebra $R$ decomposes into the tensor product of PSH algebras with only one irreducible primitive element. Precisely, let ${\mathcal C} = \Omega \bigcap P$ denote the set of irreducible primitive elements in $R$. For any $\rho \in {\mathcal C}$ set
\[ \Omega(\rho) = \{ \omega \in \Omega \ | \  \langle \omega, \rho^{n} \rangle \not= 0 \ {\rm for \ some} \ n \geq 0 \}  \]
and
\[  R(\rho) = \oplus_{\omega \in \Omega(\rho)} \ {\mathbb Z} \cdot \omega . \]
Then $R(\rho)$ is a PSH algebra with set of irreducible elements $\Omega(\rho)$, whose unique irreducible primitive is $\rho$ and 
\[ R = \otimes_{\rho \in {\mathcal C}} \ R(\rho) . \]
\end{em}
\end{theorem}

 \section{The PSH algebra of $\{ GL_{m}{\mathbb F}_{q}, \ m \geq 0 \}$}
\begin{dummy}
\label{3.0}
\begin{em}

Let $R(G)$ denote the complex representation ring of a finite group $G$. Set $R =
 \oplus_{m \geq 0} \ R(GL_{m}{\mathbb F}_{q})$ with the interpretation that  $R_{0} \cong {\mathbb Z}$,
an isomorphism which gives both a choice of unit and counit for $R$. 

Let $U_{k,m-k} \subset GL_{m}{\mathbb F}_{q}$
denote the subgroup of matrices of the form
\[  X =  \left(  
\begin{array}{cc}
I_{k} & W \\
\\
0 & I_{m-k} 
\end{array} \right) \]
where $W$ is an $k \times (m-k)$ matrix. Let $P_{k, m-k}$ denote the parabolic subgroup of 
$GL_{m}{\mathbb F}_{q}$ given by matrices obtained by replacing the identity matrices $I_{k}$ and $I_{m-k}$ in the condition for membership of $U_{k,m-k}$ by matrices from $GL_{k}{\mathbb F}_{q}$ and $GL_{m-k}{\mathbb F}_{q}$ respectively. Hence there is a group extension of the form
\[      U_{k, m-k} \longrightarrow  P_{k, m-k} \longrightarrow GL_{k}{\mathbb F}_{q} \times GL_{m-k}{\mathbb F}_{q}    . \]
If $V$ is a complex representation of $GL_{m}{\mathbb F}_{q}$ then the fixed points $V^{U_{k,m-k}}$ is a representation of $GL_{k}{\mathbb F}_{q} \times GL_{m-k}{\mathbb F}_{q} $ which gives the $(k, m-k)$ component of 
\[ m^{*} : R \longrightarrow R \otimes R . \]
Given a representation $W$ of $GL_{k}{\mathbb F}_{q} \times GL_{m-k}{\mathbb F}_{q} $ so that $W \in R_{k} \otimes R_{m-k}$ we may form
\[      {\rm Ind}_{P_{k, m-k}}^{GL_{m}{\mathbb F}_{q}}( {\rm Inf}_{GL_{k}{\mathbb F}_{q} \times GL_{m-k}{\mathbb F}_{q} }^{P_{k, m-k}}(W))  \]
which gives the $(k, m-k)$ component of 
\[  m : R \otimes R \longrightarrow  R . \]

We choose a basis for $R_{m}$ to be the irreducible representations of $GL_{m}{\mathbb F}_{q}$ so that $R^{+}$ consists of the classes of representations (rather than virtual ones). Therefore it is clear that $m, m^{*}, \epsilon, \epsilon^{*}$ satisfy positivity. The inner product on $R$ is given by the Schur inner product so that for two representations $V, W$ of $GL_{m}{\mathbb F}_{q}$ we have
\[ \langle V , W \rangle = {\rm dim}_{{\mathbb C}}({\rm Hom}_{GL_{m}{\mathbb F}_{q}}(V, W)) \]
and for $m \not= n$ $R_{n}$ is orthogonal to $R_{m}$. As is well-known, with these choice of inner product, the basis of irreducible representations for $R$ is an orthonormal basis.

The irreducible primitive elements are represented by irreducible complex representations of $GL_{m}{\mathbb F}_{q}$ which have no non-zero fixed vector for any of the subgroups $U_{k, m-k}$. These representations are usually called cuspidal.

In the remainder of this section we shall verify that $R$ is a PSH algebra, as is shown in 
(\cite{AVZ81}  Chapter III). I believe, in different terminology, this structural result was known to Sandy Green at the time of writing \cite{JAG55} and to his research supervisor Phillip Hall.
\end{em}
\end{dummy}
\begin{theorem}{(Self-adjoint)}
\label{3.1}
\begin{em}

 If $X, Y, Z$ are complex representations of $GL_{m}{\mathbb F}_{q}, GL_{n}{\mathbb F}_{q},  GL_{m+n}{\mathbb F}_{q} $ respectively then
 \[    \langle m(X \otimes Y), Z \rangle =  \langle X \otimes Y, m^{*}(Z ) \rangle .\]
 
Also $\epsilon$ and $\epsilon^{*}$ are mutually adjoint.
\end{em}
\end{theorem}
\vspace{4pt}
\newpage

{\bf Proof:}

This follows from Frobenius reciprocity (\cite{Sn94} Theorem 1.2.39) because the Schur inner product is given by
\[  \begin{array}{ll} 
 \langle m(X \otimes Y), Z \rangle & = 
{\rm dim}_{{\mathbb C}}({\rm Hom}_{GL_{m+n}{\mathbb F}_{q} }( m(X \otimes Y), Z )) \\
\\
&=  {\rm dim}_{{\mathbb C}}({\rm Hom}_{P_{m, n}}
 {\rm Inf}_{GL_{m}{\mathbb F}_{q} \times GL_{n}{\mathbb F}_{q} }^{P_{m, n}}(X \otimes Y), Z))   \\
 \\
 &=  {\rm dim}_{{\mathbb C}}({\rm Hom}_{P_{m, n}}
 {\rm Inf}_{GL_{m}{\mathbb F}_{q} \times GL_{n}{\mathbb F}_{q} }^{P_{m, n}}(X \otimes Y), Z^{U_{m,n}})) .
\end{array}    \]
The adjointness of $\epsilon$ and $\epsilon^{*}$ is obvious.
$\Box$

\begin{proposition}{(Associativity and coassociativity)}
\label{3.2}
\begin{em}

The coproduct $m^{*}$ is coassociative and the product $m$ is associative.
\end{em}
\end{proposition}
\vspace{4pt}
{\bf Proof:}

Clearly $m^{*}$ is coassociative because taking fixed-points $GL_{a}{\mathbb F}_{q} \times GL_{b}{\mathbb F}_{q} \times GL_{c}{\mathbb F}_{q}$ of a $GL_{a+b+c}{\mathbb F}_{q}$ representation is clearly associative. It follows from Theorem \ref{3.1} that $m$ is associative, since the Schur inner product is non-singular.  $\Box$

\begin{theorem}{(Hopf condition)}
\label{3.3}
\begin{em}

The homomorphism $m^{*}$ is an algebra homomorphism with respect to $m$. The homomorphism $m$ is a coalgebra homomorphism with respect to $m^{*}$. 
\end{em}
\end{theorem}

Obviously the coalgebra homomorphism assertion follows from the algebra homomorphism assertion by the adjointness property of Theorem \ref{3.1}.

The discussion which follows will establish Theorem \ref{3.3}. It is rather delicate and involved so I am going to give it in full detail (following (\cite{AVZ81}  p.167 and p.173 with minor changes). For notational convenience I shall write $G_{n} = GL_{n}{\mathbb F}_{q}$ for the duration of this discussion.

Recall that we are attempting to show that for each $(\alpha, m- \alpha)$ and $(a, m-a)$ that the $R(G_{a}) \otimes R(G_{m-a})$-component of $m^{*} \cdot m$
\[  R(G_{\alpha}) \otimes R(G_{m- \alpha}) \stackrel{m}{\longrightarrow} R(G_{m}) \stackrel{m^{*}}{\longrightarrow} R(G_{a}) \otimes R(G_{m - a}) \]  
is equal to the $R(G_{a}) \otimes R(G_{m-a})$-component 
\[   R(G_{\alpha}) \otimes R(G_{m- \alpha})   \stackrel{m^{*} \otimes m^{*} }{\longrightarrow}  
R \otimes R \otimes R \otimes R  \stackrel{1 \otimes T \otimes 1}{\longrightarrow}  R \otimes R \otimes R \otimes R   \stackrel{m \otimes m}{\longrightarrow}  R \otimes R.  \]

Let $Z$ be a complex representation of $G_{m}$ then the $(a, m-a)$-component of $m^{*}(Z)$ is given by
\[ Z^{U_{a, m-a}} \in R(G_{a} \times G_{m-a}) \cong R(G_{a}) \otimes R(G_{m-a})  \]
with the group action given by the induced $P_{a, m-a}/U_{a,m-a}$-action.

If $Z = m(X \otimes Y)$ with $X,Y$ representations of $G_{\alpha}, G_{m - \alpha}$ respectively
then 
\[    Z =  m(X \otimes Y) =  {\rm Ind}_{P_{\alpha, m-\alpha}}^{G_{m}}( {\rm Inf}_{G_{\alpha} \times G_{m-\alpha} }^{P_{\alpha, m-\alpha}}(X \otimes Y)) . \]
Therefore we must study the restriction
\[ {\rm Res}_{P_{a, m-a}}^{G_{m}}(   {\rm Ind}_{P_{\alpha, m-\alpha}}^{G_{m}}( {\rm Inf}_{G_{\alpha} \times G_{m-\alpha} }^{P_{\alpha, m-\alpha}}(X \otimes Y)))  \]
by means of the Double Coset Formula (\cite{Sn94} Theorem 1.2.40; see also Chapter 7, \S1).
Explicitly the Double Coset Formula in this case gives
\[   \sum_{g \in P_{a,m-a} \backslash G_{m}  /  P_{\alpha, m - \alpha} } \   {\rm Ind}_{P_{a, m-a} \cap  gP_{\alpha, m-\alpha}g^{-1}}^{P_{a, m-a}}((g^{-1})^{*} {\rm Inf}_{G_{\alpha} \times G_{m-\alpha} }^{P_{\alpha, m-\alpha}}(X \otimes Y)   )   \]
where the $(g^{-1})^{*}$-action is given by $(ghg^{-1})(w) = hw$. 

The Double Coset Formula isomorphism (downwards) is given by
\[  z \otimes_{ P_{\alpha, m - \alpha} } w  \mapsto  j \otimes_{ P_{a, m-a}  \cap  
gP_{\alpha, m-\alpha}g^{-1}}  h w \]
where $z = jgh$ with $j \in P_{a, m-a}, h \in P_{\alpha, m - \alpha} $ with inverse (upwards) given by
\[   j \otimes_{ P_{a, m-a}  \cap  
gP_{\alpha, m-\alpha}g^{-1}}   w   \mapsto   jg \otimes_{ P_{\alpha, m - \alpha} } w.  \]

Next let $\Sigma_{m} \subset GL_{m}{\mathbb F}_{q}$ denote the symmetric group on $m$ letters embedded as the subgroup of permutation matrices (i.e. precisely one non-zero entry on each row and column which is equal to $1$). 

The following result is proved in (\cite{AVZ81} p.173; see also \cite{NB68} Chapter IV, \S2)
\begin{theorem}{(Bruhat Decomposition)}
\label{3.4}
\begin{em}

The inclusion of  $\Sigma_{m}$ into $GL_{m}{\mathbb F}_{q}$ induces a bijection
\[    \Sigma_{a} \times  \Sigma_{m-a} \backslash \Sigma_{m}  /   \Sigma_{\alpha} \times  \Sigma_{ m - \alpha}    \stackrel{\cong}{\longrightarrow}   P_{a,m-a} \backslash G_{m}  /  P_{\alpha, m - \alpha}  \]
\end{em}
\end{theorem}

Now we shall construct a convenient set of double coset representations for the left-hand side of Theorem \ref{3.4}. 

Consider the double cosets
\[\Sigma_{a} \times \Sigma_{m-a} \backslash \Sigma_{m} / \Sigma_{\alpha} \times \Sigma_{m- \alpha} .\]
On page 171 of \cite{AVZ81} one finds the assertion that the double cosets in the title of this section are in bijection with the matrices of non-negative integers
\[   \left( \begin{array}{cc}
k_{1,1} & k_{1,2} \\
\\
k_{2,1} & k_{2,2} \\
\end{array} \right)   \]
which satisfy
\[  k_{1,1} + k_{1,2} = \alpha, \  k_{2,1} + k_{2,2} = m - \alpha, \  k_{1,1} + k_{2,1} = a, \ 
k_{1,2} + k_{2,2} = m-a  . \]

Let $w \in \Sigma_{m}$ be a permutation of $\{ 1, \ldots , m \}$. Set $I_{1} = \{ 1, 2, \ldots , a \}$, $I_{2} =  \{ a+1, a+2, \ldots , m \}$, $J_{1} =  \{ 1, 2, \ldots , \alpha \}$
and $J_{2} =  \{ \alpha + 1, \alpha + 2, \ldots , m \}$. Therefore if $g \in \Sigma_{a} \times \Sigma_{m-a}$ and $g' \in \Sigma_{\alpha} \times \Sigma_{m- \alpha}$ we have, for $t=1,2$ and $v=1,2$,
\[  gwg'(J_{t}) \bigcap I_{v} =  gw(J_{t}) \bigcap I_{v} = g( w(J_{t}) \bigcap  g^{-1}(I_{v}))  = 
g( w(J_{t}) \bigcap  I_{v}). \]
Therefore if we set
\[ k_{t,v} =  \#( w(J_{t}) \bigcap  I_{v})  \]
we have a well-defined map of sets from the double cosets to the $2 \times 2$ matrices of the form described above because
\[   k_{1,v} + k_{2,v} =   \#(I_{v}) = \left\{  \begin{array}{ll}
a & {\rm if} \ v=1 , \\
\\
m-a & {\rm if} \ v=2 
\end{array}  \right. \]
and
\[   k_{t, 1} + k_{t, 2,} =   \#(J_{t}) = \left\{  \begin{array}{ll}
\alpha & {\rm if} \ t=1 , \\
\\
m- \alpha & {\rm if} \ t=2 .
\end{array}  \right. \]
Next we consider the passage from the matrix of $k_{i,j}$'s to a double coset. Write
\[ J_{1} = J(k_{1,1}) \bigcup J(k_{1,2}) , \ J_{2} = J(k_{2,1}) \bigcup J(k_{2,2}) \]
where $J(k_{1,1}) = \{1, \ldots , k_{1,1} \}$ and $J(k_{2,1}) = \{ \alpha +1, \ldots , \alpha + k_{2,1} \}$.
Similarly write
\[ I_{1} = I(k_{1,1}) \bigcup I(k_{2,1}) , \ I_{2} = I(k_{1,2}) \bigcup I(k_{2,2}) \]
where $I(k_{1,1}) = \{1, \ldots , k_{1,1} \}$ and $I(k_{1,2}) = \{ a +1, \ldots , a + k_{1,2} \}$. Since the orders of $ I(k_{i,j}) $ and $ J(k_{i,j}) $ are both equal to $k_{i,j}$ there is a permutation, denoted by $w(k_{*,*})$ which sends $J_{1} \bigcup J_{2}$ to $I_{1} \bigcup I_{2}$ by the identity on $J(k_{1,1}) = I(k_{1,1})$ and $J(k_{2,2}) = I(k_{2,2})$ and interchanges $J(k_{1,2}), J(k_{2,1})$ with $I(k_{1,2}), I(k_{2,1})$ in an order-preserving manner.

Given the permutation $w(k_{*,*})$ we have
\[  \begin{array}{l}
\#( w(k_{*,*})(J_{1}) \bigcap  I_{1})  =   \#( J(k_{1,1}) \bigcap  I(k_{1,1}) ) = k_{1,1} , \\
\\
\#( w(k_{*,*})(J_{1}) \bigcap  I_{2})  =   \#(  w(k_{*,*})(J(k_{2,1})) \bigcap  I(k_{1,2}) ) = k_{1,2} , \\
\\
\#( w(k_{*,*})(J_{2}) \bigcap  I_{1})  =   \#(  w(k_{*,*})(J(k_{1,2})) \bigcap  I(k_{2,1}) ) = k_{2,1} , \\
\\
\#( w(k_{*,*})(J_{2}) \bigcap  I_{2})  =   \#( J(k_{2,2}) \bigcap  I(k_{2,2}) ) = k_{2,2} 
\end{array}  \]
so that the map $k_{*,*} \mapsto  \Sigma_{a} \times \Sigma_{m-a}  w(k_{*,*})  \Sigma_{\alpha} \times \Sigma_{m- \alpha}$ is a split injection. In addition it is straightforward to verify that any permutation whose $k_{*,*}$-matrix equals that of $w(k_{*,*})$ belongs to the same double coset as $w(k_{*,*})$. Hence the map is a bijection.

For example when $a=3, \alpha = 4, k_{11} = 1 = k_{22}, k_{21}=2, k_{12}=3$
\[   w(k_{*,*})^{-1} =   \left(  \begin{array}{ccccccc}
1 & 0 & 0 & 0&0&0 &0 \\
 0 &0& 0 & 1&0&0&0 \\
 0 &0&0&0&1& 0 &0 \\
 0 & 0 & 0 & 0&0&1&0 \\
0&1&0&0&0&0&0  \\
0&0&1&0&0& 0&0  \\
0&0&0&0&0&0&1 \\
\end{array}  \right)  \]
\[    w(k_{*,*}) =   \left(  \begin{array}{ccccccc}
1 & 0 & 0 & 0 &0 &0 & 0 \\
 0 & 0& 0 & 0&1&0&0 \\
0 &0&0&0&0& 1 &0 \\
 0 & 1 & 0 & 0&0&0&0 \\
0&0&1&0&0&0&0 \\
0&0&0&1&0& 0&0 \\
0&0&0&0&0&0&1 \\
\end{array}  \right) .  \]

This permutation arises in another way as a permutation of the basis elements of tensor products of four vector spaces. Let
\[  V_{1} = {\mathbb F}_{q}^{k_{11}} \oplus  {\mathbb F}_{q}^{k_{12}} \oplus  {\mathbb F}_{q}^{k_{21}} \oplus  {\mathbb F}_{q}^{k_{22}}  \ {\rm and}   \    V_{2} = {\mathbb F}_{q}^{k_{11}} \oplus  {\mathbb F}_{q}^{k_{21}} \oplus  {\mathbb F}_{q}^{k_{12}} \oplus  {\mathbb F}_{q}^{k_{22}} . \]
We have the linear map 
\[  1 \oplus T(k_{*,*}) \oplus 1 : V_{1} \longrightarrow  V_{2} \]
which interchanges the order of the two central direct sum factors. The basis for 
$V_{1}$ is made in the usual manner from ordered bases $\{ e_{1}, \ldots \ e_{k_{11}} \}$,  $\{ e_{k_{11}+1}, \ldots \ e_{k_{11}+ k_{12}} \}$, $\{ e_{k_{11}+ k_{12}+1}, \ldots \ e_{k_{11}+ k_{12}+ k_{21}} \}$
\linebreak
and $\{ e_{k_{11}+ k_{12}+ k_{21}+1}, \ldots \ e_{m} \}$ of ${\mathbb F}_{q}^{k_{11}}$, ${\mathbb F}_{q}^{k_{12}}$,
$  {\mathbb F}_{q}^{k_{21}}$ and $  {\mathbb F}_{q}^{k_{22}}$ respectively. Similarly the basis for 
$V_{2}$ is made in the usual manner from ordered bases $\{ v_{1}, \ldots \ v_{k_{11}} \}$,  $\{ v_{k_{11}+1}, \ldots \ v_{k_{11}+ k_{21}} \}$, $\{ v_{k_{11}+ k_{21}+1}, \ldots \ v_{k_{11}+ k_{21}+ k_{12}} \}$
\linebreak
and $\{ v_{k_{11}+ k_{21}+ k_{12}+1}, \ldots \ v_{m} \}$ of ${\mathbb F}_{q}^{k_{11}}$, ${\mathbb F}_{q}^{k_{21}}$,
$  {\mathbb F}_{q}^{k_{12}}$ and $  {\mathbb F}_{q}^{k_{22}}$ respectively.

The linear map $1 \oplus T(k_{*,*}) \oplus 1$ sends the ordered set  $\{ e_{1}, \ldots , e_{m} \}$ to the order set $\{ v_{1}, \ldots , v_{m} \}$ by $e_{j} \mapsto v_{ w(k_{*,*})(j)}$.

Clearly
\[   w(k_{*,*})G_{k_{11}} \times  G_{k_{12}} \times G_{k_{21}} \times G_{k_{22}}w(k_{*,*})^{-1} = 
G_{k_{11}} \times  G_{k_{21}} \times G_{k_{12}} \times G_{k_{22}}   \]
from which is it easy to see that
\[   w(k_{*,*})G_{\alpha} \times  G_{m - \alpha}w(k_{*,*})^{-1}  \bigcap G_{a} \times G_{m-a} = 
G_{k_{11}} \times  G_{k_{21}} \times G_{k_{12}} \times G_{k_{22}}   . \]

So far we have shown that the $(a, m-a)$-component of $m^{*}(m(X \otimes Y))$ is the sum of terms, one for each $w(k_{*,*})$, given by the induced $G_{a} \times G_{m-a}$-action on
\[      {\rm Ind}_{P_{a, m-a} \cap  w(k_{*,*}P_{\alpha, m-\alpha}w(k_{*,*})^{-1}}^{P_{a, m-a}}(( w(k_{*,*})^{-1})^{*} {\rm Inf}_{G_{\alpha} \times G_{m-\alpha} }^{P_{\alpha, m-\alpha}}(X \otimes Y)   )  .  \]

On the other hand, for each $w(k_{*,*})$ there is a $(a, m-a)$-component of the other composition we are studying given by
\[  \begin{array}{l}
R(G_{\alpha} \times G_{m - \alpha}) 
\stackrel{(-)^{U_{k_{11}, k_{12}} \times U_{k_{21},k_{22} }}}{\longrightarrow}
R(G_{k_{11}} \times  G_{k_{12}} \times G_{k_{21}} \times G_{k_{22}} ) \\
\\
\hspace{60pt} \stackrel{1 \otimes T(k_{*,*}) \otimes 1}{\longrightarrow} 
R(G_{k_{11}} \times  G_{k_{21}} \times G_{k_{12}} \times G_{k_{22}} ) \\
\\
\hspace{80pt} \stackrel{{\rm IndInf} \times {\rm IndInf}}{\longrightarrow}   R(G_{a} \times G_{m - a}) .
\end{array}  \]
Composing this second route with the split surjection
\[  R(G_{a} \times G_{m - a}) \stackrel{{\rm Inf}}{\longrightarrow} R(P_{a,m-a}) \]
is equal to the composition
\[  \begin{array}{l}
R(G_{\alpha} \times G_{m - \alpha}) 
\stackrel{(-)^{U_{k_{11}, k_{12}} \times U_{k_{21},k_{22} }}}{\longrightarrow}
R(G_{k_{11}} \times  G_{k_{12}} \times G_{k_{21}} \times G_{k_{22}} ) \\
\\
\hspace{20pt} \stackrel{1 \otimes T(k_{*,*}) \otimes 1}{\longrightarrow} 
R(G_{k_{11}} \times  G_{k_{21}} \times G_{k_{12}} \times G_{k_{22}} ) \stackrel{{\rm Inf}}{\longrightarrow} \\
\\
\hspace{40pt}  R(P_{k_{11}, k_{21}, k_{12}, k_{22}})  \stackrel{{\rm Ind}}{\longrightarrow}   
R(P_{a, m - a}) 
\end{array}  \]
because the kernels of the quotient maps $ P_{k_{11}, k_{21}, k_{12}, k_{22}} \longrightarrow   P_{k_{11}, k_{21}} $ and 
\linebreak
$  P_{a, m-a} \longrightarrow  G_{a} \times G_{m-a} $ are both equal to $U_{a,m-a}$.

This composition takes the $U_{k_{11}, k_{12}} \times  U_{k_{21}, k_{22}}$-fixed points of $X \otimes  Y$ with the $G_{k_{11}} \times  G_{k_{12}} \times  G_{k_{21}} \times  G_{k_{22}} $-action and then conjugates it by $w(k_{*,*})$. Alternatively it takes the $w(k_{*,*})U_{k_{11}, k_{12}} \times  U_{k_{21}, k_{22}}w(k_{*,*})^{-1}$-fixed points of $(w(k_{*,*})^{-1})^{*}(X \otimes  Y)$ with the $G_{k_{11}} \times  G_{k_{21}} \times  G_{k_{12}} \times  G_{k_{22}} $-action. Now 
\[  w(k_{*,*})U_{k_{11}, k_{12}} \times  U_{k_{21}, k_{22}}w(k_{*,*})^{-1} \subset U_{a, m-a}.    \]
For example, in the small example given in the last section of this chapter, $U_{k_{11}, k_{12}} \times  U_{k_{21}, k_{22}}$ consists of matrices of the form
\[ D = \left(  \begin{array}{ccccccc}
1&a_{12}&a_{13}& a_{14} &0&0&0\\
0&1&0&0&0&0&0 \\
0&0&1&0&0&0&0\\
0&0&0&1&0&0&0\\
0&0&0&0&1&0&a_{57} \\
0&0&0&0&0&1&a_{67}\\
0&0&0&0&0&0&1\\
\end{array} \right)  \]
so that $w(k_{*,*})U_{k_{11}, k_{12}} \times  U_{k_{21}, k_{22}}w(k_{*,*})^{-1}$ consists of matrices
\[  w(k_{*,*})Dw(k_{*,*})^{-1} = \left(  \begin{array}{ccccccc}
1&0&0&a_{12}&a_{13}& a_{14} &0\\
0&1&0&0&0&0&a_{57} \\
0&0&1&0&0&0&a_{67}\\
0&0&0&1&0&0&0\\
0&0&0&0&1&0&0 \\
0&0&0&0&0&1&0\\
0&0&0&0&0&0&1\\
\end{array} \right)  .  \]
Since $ w(k_{*,*})Dw(k_{*,*})^{-1}$'s act trivially we may inflate the representation to 
\newline
$P_{k_{11}, k_{21}, k_{12}, k_{22}}$ (i.e. extending the action trivially on $U_{k_{11}, k_{21}, k_{12}, k_{22}}$) and then induce up to a representation of $P_{a, m-a}$.

Now let us describe the isomorphism between the result of sending $X \otimes Y$ via the second route and the $U_{a, m-a}$-fixed subspace of 
\[      {\rm Ind}_{P_{a, m-a} \cap  w(k_{*,*}P_{\alpha, m-\alpha}w(k_{*,*})^{-1}}^{P_{a, m-a}}(( w(k_{*,*})^{-1})^{*} {\rm Inf}_{G_{\alpha} \times G_{m-\alpha} }^{P_{\alpha, m-\alpha}}(X \otimes Y)   )  .  \]

There are inclusions
\[ w(k_{*,*})P_{\alpha, m- \alpha}w(k_{*,*})^{-1} \bigcap P_{a,m-a} \subset  P_{k_{11}, k_{21}, k_{12}, k_{22}} \subset P_{a, m-a} .\]
For example, $w(k_{*,*})P_{\alpha, m - \alpha}w(k_{*,*})^{-1}  \bigcap P_{3,4} $, in the small example of the Appendix, consists of the matrices of the form
\[   E' = \left( \begin{array}{ccccccc}
a_{11} &  a_{15}&  a_{16}&  a_{12}&  a_{13}&  a_{14}&  a_{17}   \\
 0&  a_{55}&  a_{56}&  0&  0&  0&  a_{57}  \\
0&  a_{65}&  a_{66}&  0&  0&  0&  a_{67}  \\
 0&  0&  0& a_{22}&  a_{23}&  a_{24}&  a_{27}  \\
 0&  0&  0&  a_{32}&  a_{33}&  a_{34}&  a_{37} \\
0 &  0&  0&  a_{42}&  a_{43}&  a_{44}&  a_{47} \\
0&  0&  0&  0&  0&  0&  a_{77}  \\
\end{array}  \right)\]
and, as we noted above, $ w(k_{*,*})Dw(k_{*,*})^{-1}$ consists of the matrices
\[  \{ (b_{ij}) \in w(k_{*,*})P_{\alpha, m- \alpha}w(k_{*,*})^{-1} \bigcap U_{3,4} \ | \  b_{17}=0  \} . \]

There is a bijection of cosets
\[  \begin{array}{l}
P_{k_{11}, k_{21}, k_{12}, k_{22}}/w(k_{*,*})P_{\alpha, m- \alpha})w(k_{*,*})^{-1} \bigcap P_{a,m-a} \\
\\
\cong U_{a, m-a}/  w(k_{*,*})P_{\alpha, m- \alpha})w(k_{*,*})^{-1} \bigcap U_{3,4} .
\end{array}  \]
Therefore we may take the coset representations $X_{\alpha}$ to lie in the abelian group $U_{a, m-a}$.
In the small example the $X_{\alpha}$'s may be taken to be of the form
\[   X_{\alpha} = \left( \begin{array}{ccccccc}
1 &  0&0&0&0&0&a  \\
 0&  1&  0&  b&  c&  d&  0  \\
0&  0&  1&  e&  f&  g& 0  \\
 0&  0&  0& 1&  0& 0&  0  \\
 0&  0&  0&  0&  1&  0&  0 \\
0 &  0&  0& 0& 0&  1&  0 \\
0&  0&  0&  0&  0&  0&  1  \\
\end{array}  \right)  .    \]
The isomorphism from the image of $X \otimes Y$ via the second route to the $U_{a, m-a}$-fixed subspace of 
\[      {\rm Ind}_{P_{a, m-a} \cap  w(k_{*,*}P_{\alpha, m-\alpha}w(k_{*,*})^{-1}}^{P_{a, m-a}}(( w(k_{*,*})^{-1})^{*} {\rm Inf}_{G_{\alpha} \times G_{m-\alpha} }^{P_{\alpha, m-\alpha}}(X \otimes Y)   )    \]
is given by
 \[  g \otimes_{P_{k_{11}, k_{21}, k_{12}, k_{22}}} v \mapsto 
  \sum_{X_{\alpha}} \  g X_{\alpha} \otimes_{w(k_{*,*})P_{\alpha, m- \alpha})w(k_{*,*})^{-1} \bigcap P_{a,m-a} } v .  \]
  
This concludes the proof of Theorem \ref{3.3}. The remainder of the Hopf condition is given by the following result, which is proved in a similar manner to Theorem \ref{3.3} (see \cite{AVZ81} p.175).
\begin{theorem}{$_{}$}
\label{3.5}
\begin{em}

In the notation of \S1, $\epsilon^{*}$ is a ring homomorphism restricting to an isomorphism on $R_{0}$.
\end{em}
\end{theorem} 
\section{The PSH algebra of a wreath product}

Let $H$ be a finite group and  consider the wreath-product $\Sigma_{n} \int H$. This is the semi-direct product of the symmetric group $\Sigma_{n}$ acting on the $n$-fold product $H^{n} = H \times \ldots \times H$ on the left by permuting the factors.
An arbitrary element of this group is $(\sigma , \alpha_{1}, \ldots , \alpha_{n})$ where $\sigma \in \Sigma_{n}$ and $\alpha_{i} \in H$. The product is given by
\[ (\sigma , \alpha_{1}, \ldots , \alpha_{n}) \cdot (\sigma' , \alpha'_{1}, \ldots , \alpha'_{n}) =
( \sigma \sigma' , \alpha_{1} \alpha'_{\sigma(1)} , \alpha_{2} \alpha'_{\sigma(2)} , \ldots , \alpha_{n} \alpha'_{\sigma(n)}).  \]

We pause to check associativity and inverses, since we shall need the explicit formulae later.
\[ \begin{array}{l}
((\sigma , \alpha_{1}, \ldots , \alpha_{n}) \cdot (\sigma' , \alpha'_{1}, \ldots , \alpha'_{n})) \cdot 
(\sigma'', \alpha_{1}'', \ldots , \alpha_{n}'') \\
\\
= ( \sigma \sigma' , \alpha_{1} \alpha'_{\sigma(1)} , \alpha_{2} \alpha'_{\sigma(2)} , \ldots , \alpha_{n} \alpha'_{\sigma(n)}) \cdot (\sigma'', \alpha_{1}'', \ldots , \alpha_{n}'') \\
\\
= (  \sigma \sigma' \sigma'',  \alpha_{1} \alpha'_{\sigma(1)}   \alpha_{ \sigma \sigma' (1)}'' , \alpha_{2} \alpha'_{\sigma(2)}  \alpha_{ \sigma \sigma' (2)}'' , \ldots , \alpha_{n} \alpha'_{\sigma(n)}  \alpha_{ \sigma \sigma' (n)}'' )  \\
\\
{\rm and}  \\
\\
(\sigma , \alpha_{1}, \ldots , \alpha_{n}) \cdot ((\sigma' , \alpha'_{1}, \ldots , \alpha'_{n}) \cdot 
(\sigma'', \alpha_{1}'', \ldots , \alpha_{n}'')) \\
\\
= (\sigma , \alpha_{1}, \ldots , \alpha_{n}) \cdot  (\sigma' \sigma'', \alpha'_{1} \alpha_{\sigma'(1)}'', \ldots , \alpha'_{n} \alpha_{\sigma'(n)}'') \\
\\
=  (  \sigma \sigma' \sigma'',  \alpha_{1} \alpha'_{\sigma(1)}   \alpha_{ \sigma \sigma' (1)}'' , \alpha_{2} \alpha'_{\sigma(2)}  \alpha_{ \sigma \sigma' (2)}'' , \ldots , \alpha_{n} \alpha'_{\sigma(n)}  \alpha_{ \sigma \sigma' (n)}'' ) .
\end{array} \]
Also
\[   \begin{array}{l}
(\sigma  , \alpha_{1}, \ldots   , \alpha_{n})^{-1} = (\sigma^{-1}  , (\alpha_{(\sigma)^{-1}(1)})^{-1}, \ldots , ( \alpha_{(\sigma)^{-1}(n)})^{-1}) \\
\\
{\rm since} \\
\\
(\sigma, \alpha_{1}, \ldots , \alpha_{n}) \cdot  (\sigma^{-1}  , (\alpha_{(\sigma)^{-1}(1)})^{-1}, \ldots , (\alpha_{(\sigma)^{-1}(n)})^{-1}) \\
\\
= ( \sigma \sigma^{-1}  , \alpha_{1}  (\alpha_{\sigma \sigma^{-1}(1)})^{-1}  , \ldots ,  \alpha_{n}  (\alpha_{\sigma \sigma^{-1}(n)})^{-1}  ) \\
\\
= (1, 1, \ldots , 1).
\end{array} \]

The PSH algebra associated to complex representations of these wreath products are defined in (\cite{AVZ81} Chapter Two \S7). The underlying graded group is 
\[ R = \oplus_{n \geq 0} \ R(\Sigma_{n} \int H)  . \]

The multiplication
\[ m :  R(\Sigma_{m} \int H) \otimes  R(\Sigma_{n} \int H) \longrightarrow   R(\Sigma_{m+n} \int H)  \]
is defined on representations by 
\[ m(V \otimes W) = {\rm Ind}_{\Sigma_{m} \int H) \times \Sigma_{n} \int H)}^{\Sigma_{m+n} \int H)}(V \otimes W).\]
The comultiplication 
\[ m^{*} : R(\Sigma_{n} \int H))  \longrightarrow  \oplus_{a=0}^{n} \ R(\Sigma_{a} \int H)) \otimes R(\Sigma_{n-a} \int H))  \]
is defined to have $(a, n-a)$ component given by
\[ {\rm Res}_{\Sigma_{a} \int H \times  \Sigma_{n-a} \int H}^{\Sigma_{n} \int H} :  R(\Sigma_{n} \int H)  \longrightarrow   R(\Sigma_{a} \int H) \otimes R(\Sigma_{n-a} \int H) . \]

By convention $R_{0} = {\mathbb Z}$ and the unit and counit are defined as in the 
$GL_{n}{\mathbb F}_{q}$ example of \S3. 

When $H$ is trivial this is the PSH algebra associated to the symmetric groups.

The set of irreducible primitives $\Omega$ is given by the set of representations of $H$, considered as positive elements in $R(\Sigma_{1} \int H) = R(H)$.

\section{Kondo-Gauss sums}

\begin{definition}
\label{5.1}
\begin{em}

Let $\rho : H \longrightarrow  GL_{n}{\mathbb C}$ denote a representation of a subgroup $H$ of $GL_{n}{\mathbb F}_{q}$. If $q$ is a power of the prime $p$ we have the (additive) trace map
\[    {\rm Tr}_{{\mathbb F}_{q}/{\mathbb F}_{p}} : {\mathbb F}_{q} \longrightarrow {\mathbb F}_{p}.   \]
In addition we have the matrix trace map
\[ {\rm Trace}  : GL_{n}{\mathbb F}_{q} \longrightarrow  {\mathbb F}_{q}  . \]
Define a measure map $\Psi$ on matrices $X \in  GL_{n}{\mathbb F}_{q} $ by
\[     \Psi(X) =  e^{  \frac{2 \pi \sqrt{-1}  {\rm Tr}_{{\mathbb F}_{q}/{\mathbb F}_{p}}( {\rm Trace}(X)) }{p}    }  \]
which is denoted by $e_{1}[X]$ in \cite{Kon63}. Let $\chi_{\rho}$ denote the character function of $\rho$ which assigns to $X$ the trace of the complex matrix $\rho(X)$.

Define a complex number $W_{H}(\rho)$ by the formula
\[ W_{H}(\rho) = \frac{1}{ {\rm dim}_{{\mathbb C}}(\rho) }   \  \sum_{X \in H} \    \chi_{\rho}(X) \Psi(X)  .\]

When $H = GL_{n}{\mathbb F}_{q}$ and $\rho$ is irreducible $W_{GL_{n}{\mathbb F}_{q}}(\rho) =  w(\rho)$, the Kondo-Gauss sum which is introduced and computed in \cite{Kon63}.
\end{em}
\end{definition}
\begin{theorem}{(\cite{Sn18} Appendix III Theorem 3.2)}
\label{5.2}
\begin{em}

Let $\sigma $ be a finite-dimensional representation of $H \subseteq GL_{n}{\mathbb F}_{q}$. Then for any subgroup $J$
such that $H \subseteq J  \subseteq GL_{n}{\mathbb F}_{q}$
\[   W_{H}(\sigma) = W_{J}( {\rm Ind}_{H}^{J}(\sigma)) .\]
\end{em}
\end{theorem}
\vspace{2pt}

{\bf Proof}
\vspace{2pt}

Set $\rho =  {\rm Ind}_{H}^{J}(\sigma)$. By definition
\[ \begin{array}{ll}
W_{J}(\rho) &=    \frac{|H|}{|J| \cdot  {\rm dim}_{{\mathbb C}}(\sigma) }  \   \sum_{X \in J} \  \chi_{\rho}(X) \Psi(X)   \\
\\
& =      \frac{1}{|J| \cdot  {\rm dim}_{{\mathbb C}}(\sigma) }  \   \sum_{X \in J} \  \sum_{Y \in J, \ YXY^{-1} \in H}  \   \chi_{\sigma}(YXY^{-1}) \Psi(X) 
\end{array}  \]
by the character formula for an induced representation (\cite{Sn94} Theorem 1.2.43). Consider the free action of $J$ on $J \times J$ given by $(X,Y) Z = (Z^{-1}XZ , YZ)$ for $XY,Z \in J$. The map from $J \times J$ to $J$ sending $(X,Y)$ to $YXY^{-1}$ is constant on each $J$-orbit. Therefore
\[ \begin{array}{ll}
W_{J}(\rho) 
& =      \frac{1}{|J| \cdot  {\rm dim}_{{\mathbb C}}(\sigma) }  \   \sum_{X \in J} \  \sum_{Y \in J, \ YXY^{-1} \in H}  \   \chi_{\sigma}(YXY^{-1}) \Psi(YXY^{-1})  \\
\\
& =     \frac{1}{|J| \cdot  {\rm dim}_{{\mathbb C}}(\sigma) }  \   |J|  \sum_{U  \in H}  \   \chi_{\sigma}(U) \Psi(U)  \\
\\
& =  W_{H}(\sigma) .
\end{array}  \]
\begin{lemma}{$_{}$}
\label{5.3}
\begin{em}

\[ ({\rm dim}_{\mathbb C}(\sigma_{1}) +  {\rm dim}_{\mathbb C}(\sigma_{2})) W_{H}(\sigma_{1} \oplus \sigma_{2}) = {\rm dim}_{\mathbb C}(\sigma_{1})W_{H}(\sigma_{1} ) +   {\rm dim}_{\mathbb C}(\sigma_{2})W_{H}(\sigma_{2}) .\]
\end{em}
\end{lemma}
\begin{example}{The Weil representation $r(\Theta)$ of $GL_{2}{\mathbb F}_{q}$}
\label{5.4}
\begin{em}

The Weil representation is a very ingenious construction of a $(q-1)$-dimensional irreducible complex representation of $GL_{2}{\mathbb F}_{q}$. It is constructed from scratch in (\cite{Sn94} Chapter Three).
However there is a very simple description of $r(\Theta)$ in terms of induced representations, which may be verified (for example) using the character formulae of (\cite{Sn94} Chapter Three).

There is a copy of ${\mathbb F}_{q^{2}}^{*}$, unique up to conjugation, embedded in $GL_{2}{\mathbb F}_{q}$. For example, if $\sigma$ is a non-square in ${\mathbb F}_{q}^{*}$ then sending $a + b \sqrt{\sigma}$ to
\[   \left(  \begin{array}{cc}
a & b \sigma \\
\\
b & a  
\end{array}  \right)   \]
gives such an embedding. Let $F$ denote the generator of ${\rm Gal}({\mathbb F}_{q^{2}}/{\mathbb F}_{q})$ and suppose that $\Theta :  {\mathbb F}_{q^{2}}^{*} \longrightarrow {\mathbb C}^{*}$ is a 
character such that $F(\Theta) \not= \Theta$. 

Let $H$ denote the ``top line'' subgroup consisting of of matrices of the form
\[   \left(  \begin{array}{cc}
a & b  \\
\\
0 & 1  
\end{array}  \right)   \] 
so that $H \cong {\mathbb F}_{q}^{*} \times {\mathbb F}_{q}$ by sending the above matrix to $(a, b/a)$.
Therefore we have a character on $H$ given by 
\[  (\Theta \otimes \Psi) \left(  \begin{array}{cc}
a & b  \\
\\
0 & 1  
\end{array}  \right)  = \Theta(a) \Psi(b/a) .  \]
Here $\Psi$ is the additive measure defined in Definition \ref{5.1} on $n \times n$ matrices in the case $n=1$.

 There is a (split) short exact sequence of complex 
 $GL_{2}{\mathbb F}_{q}$-representations\footnote{To establish this result I first used the fact that (\cite{Sn94} Chapter Three) provides an easy description of the right-hand map together with a complicated argument to show that the left-hand representation was inside the kernel. Then, smugly pleased with the discovery, I check it using the character values of (\cite{Sn94} Chapter Three) only to find the same result appears in (\cite{BH06} p.47)!}
\[   0   \longrightarrow  {\rm Ind}_{{\mathbb F}_{q^{2}}^{*}}^{GL_{2}{\mathbb F}_{q}}( \Theta ) 
\longrightarrow  {\rm Ind}_{H}^{GL_{2}{\mathbb F}_{q}}( \Theta \otimes \Psi)  \longrightarrow r(\Theta)  \longrightarrow  0.\]

By Theorem \ref{5.2} and Lemma \ref{5.3} we have
\[   (q-1) W_{GL_{2}{\mathbb F}_{q}}(r(\Theta)) + (q^{2} - q) W_{{\mathbb F}_{q^{2}}^{*}}(\Theta) = (q^{1} - 1) W_{H}(\Theta \otimes \Psi) . \] 
However
\[  W_{H}(\Theta \otimes \Psi)  = \sum_{(a,b) \in H} \   \Theta(a)\Psi(b/a) \Psi(a + 1) = 0 \]
since the sum of the values of a non-trivial character over a finite abelian group (${\mathbb F}_{q}$ in this case) is zero. Therefore\footnote{The factor $q$ was inadvertently missed out in (\cite{Sn18} Appendix III Example 3.4).}
\[   W_{GL_{2}{\mathbb F}_{q}}(r(\Theta)) =  - q  W_{{\mathbb F}_{q^{2}}^{*}}(\Theta)  \] 
where the right side is $\frac{-q}{q^{2}-1}$ times the classical Gauss sum over a finite field.
\end{em}
\end{example}

\begin{proposition}{$_{}$}
\label{5.5}
\begin{em}

For $i= 1,2$ let $\sigma_{i} : H_{i} \longrightarrow GL_{n_{i}}{\mathbb C}$ be a representation of $H_{i} \subseteq GL_{s_{i}}{\mathbb F}_{q}$. Then we have a representation of $H_{1} \times H_{2}$ (embedded into $GL_{s_{1} + s_{2}}$ by direct sum of matrices) given by the tensor product $\sigma_{1} \otimes \sigma_{2}$ and
\[   W_{H_{1} \times H_{2}}(\sigma_{1} \otimes \sigma_{2}) =  W_{H_{1} }(\sigma_{1})  W_{ H_{2}}(\sigma_{2}) . \]
\end{em}
\end{proposition}
\vspace{2pt}

{\bf Proof}
\vspace{2pt}

We have
\[ \begin{array}{l}
 W_{H_{1} \times H_{2}}(\sigma_{1} \otimes \sigma_{2}) \\
 \\
 =  \frac{1}{s_{1} \cdot s_{2}} \ \sum_{X_{1} \oplus X_{2} \in H_{1} \times H_{2} } \  
\chi_{\sigma_{1} \otimes \sigma_{2}}(X_{1} \oplus X_{2}) \Psi(X_{1} \oplus X_{2}) \\
\\
=  \frac{1}{s_{1} \cdot s_{2}} \ \sum_{X_{1} \oplus X_{2} \in H_{1} \times H_{2} } \  
\chi_{\sigma_{1}}(X_{1}) \chi_{ \sigma_{2}}(X_{2})  \Psi(X_{1})  \Psi(X_{2})  \\
\\
=  W_{H_{1} }(\sigma_{1})  W_{ H_{2}}(\sigma_{2}) . 
\end{array} \]
$\Box$
\begin{remark}
\label{5.6}
\begin{em}

In  \cite{Kon63} Kondo gives formulae for his Gauss sums on the irreducible complex representations of $GL_{n}{\mathbb F}_{q}$. The calculations of \cite{Kon63} do not use the function $W_{H}(\rho)$ but stick to the case of an irreducible $\rho$ and $H = GL_{n}{\mathbb F}_{q}$. The greater freedom and generality of  $W_{H}(\rho)$ should make the calculations much simpler. In (\cite{Sn18} Appendix III) I made some calculations towards obtaining Kondo's formulae by this approach.
\end{em}
\end{remark}

 \section{Symmetric groups or $GL_{n}$ of the trivial field}

The Kondo-Gauss sums of representations of $GL_{n}{\mathbb F}_{q}$ are integrals using the measure given by the additive character $\Psi$. In this section I shall establish the analogue of the Kondo-Gauss sum results for the PSH associated to symmetric groups. In the next section I shall extend these results to wreath products. 

First we need the correct measure, meaning one obeying the propoerties proved for 
$GL_{n}{\mathbb F}_{q}$ in \S5. On pp.307-8 of \cite{Sn18} we find the Kondo-Gauss sum for the symmetric group (aka the general linear group of the ``field of one element''), which uses an disjoint-multiplicative character on the symmetric group. I learnt the  measure which is used from Francesco Mezzadri\footnote{A few years ago Francesco Mezzadri, who is my son-in-law, showed me his conjectural formula for the ``Kondo-Gauss sum'' for symmetric groups and asked me whether I could prove it. The contents of this section contain, inter alia, the proof which I came up with - based on Gordon James' book \cite{GDJ}. Alas, in a bout of stereotypical absentmindedness, I forgot to show Francesco the idea. I was only reminded of the topic when I proved the relation of the formula to the PSH algebra (\cite{Sn18} Appendix III Proosition 4.3) by which time Francesco and his colleagues had found the proof given in \cite{KMS2011}. }.

We shall require some of the fundamental facts about irreducible representations of the symmetric groups, taken from \cite{GDJ}. These are adsolutely irreducible so, although for the numerical invariants we are interested in complex representation we may as well discuss ${\mathbb Q}$-representations at this point.

\begin{definition}
\label{6.1}
\begin{em}

If $\lambda = (\lambda_{1}, \ldots , \lambda_{r}) $ is a partition of the integer $n$ by strictly positive integers, $\{ \lambda_{i} , 1 \leq i \leq r \}$ with $\lambda_{i} \geq \lambda_{i+1}$, then the diagram of $\lambda$ is denoted by $[ \lambda ]$ and consists of the integer pairs 
\[   \{   (i,j) \  |  \  1 \leq i, 1 \leq j \leq \lambda_{i} \}  .     \]
Each point $(i,j)$ is called a node of $[ \lambda ]$. For example the diagram of $(4,2,2,1)$ is depicted\footnote{The reader will notice that, in the depiction of {\em diagrams}, it is the convention that $i$'s increase to the right and $j$'s increase downward. This convention is slightly weird but that's tradition for you! When we come to tableaux, which involves putting integers on the nodes/crosses those numbers can go into arbitrary positions.} as 
\[ \begin{array}{cccc}
\times & \times & \times & \times \\
\times & \times &&  \\
\times & \times &&  \\
\times &&&  \\
\end{array} \]
which may usually be written $(4, 2^{2}, 1)$ using exponent notation for brevity.

If $\lambda$ and $\mu$ are partitions of $n$ we say that $\lambda$ dominates $\mu$, written $\lambda \unrhd \mu$ if for all $j$ 
\[   \sum_{i=1}^{j} \lambda_{i} \geq   \sum_{i=1}^{j} \mu_{i}. \]
For example, we have $(6) \unrhd (5,1) \unrhd (4,2) \unrhd (3,3) \unrhd (3,1^{3}) $.

There is also the lexicographical total order on partitions of $n$ where $\lambda > \mu$ if and only if the least $j$ for which $\lambda_{j} \not= \mu_{j}$ satisfies $\lambda_{j} > \mu_{j}$.

The total order refines the partial order in the sense that $\lambda_{j} \rhd \mu_{j}$ implies that $\lambda_{j} > \mu_{j}$ but not conversely.
\end{em}
\end{definition}
\begin{definition}
\label{6.2}
\begin{em}

If $[\lambda]$ is a diagram the conjugate diagram $[\lambda']$ is obtained by interchanging the rows and columns in $[\lambda]$ and $\lambda'$ is the partition of $n$ which is conjugate to $\lambda$. Hence $4[(4,2^{2},1)']$ looks like
\[ \begin{array}{cccc}
\times & \times & \times & \times \\
\times & \times & \times &  \\
\times &  &&  \\
\times &&&  \\
\end{array} \]
The lexicographical order is reversed upon taking conjugates and also
\[ \lambda \rhd \mu \Leftrightarrow   \mu'  \rhd  \lambda'   . \]

A $\lambda$-tableau is one of the $n!$ arrays of $n$ distinct integers obtained by replacing each node of $[\lambda]$ by one of the integers in the set $\{ 1,2,\ldots , n\}$. The symmetric group acts in such a way that $\sigma \in \Sigma_{n}$ replaces each integer $i$ in a $\lambda$-tableau by $\sigma(i)$. If $t$ is a $\lambda$-tableau then I shall denote this by $\sigma(t)$.

For example, the following are two $(4,3,1)$-tableaux
\[ \begin{array}{ccccccccccc}
1 & 2 & 4 & 5 & & {\rm and} & & 4 & 5 & 7 & 3 \\
3 & 6 & 7 & &&&&2&1&8& \\
8 &  && &&&&6&&& \\
\end{array} \]
and $\sigma = (1,4,7,8,6)(2,5,3)$ sends the left-hand tableau into the right-hand one.

If we think of a $\lambda$-tableau as giving us the permutation with cycle decomposition given by writing each row as a cyclic permutation and stringing them out in order from the top - so that the left-hand one becomes $(1,2,4,5)(3,6,7)(8)$ then $\sigma$ applied to this gives the permutation conjugated by $\sigma$. That is, in the above example 
\[   \sigma \cdot (1,2,4,5)(3,6,7)(8) \cdot \sigma^{-1} = (4,5,7,3)(2,1,8)(6)  .   \]
\end{em}
\end{definition}
\begin{lemma}{(Basic Combinatorial Lemma (\cite{GDJ} \S3.7)}
\label{6.3}
\begin{em}

Let $\lambda$ and $\mu$ be partitions of $n$ and suppose that $t_{1}$ is a $\lambda$-tableau and $t_{2}$ is a $\mu$-tableau. Suppose that for each $i$ the numbers from the $i$-th row of $t_{2}$ belong to different columns of $t_{1}$ then $\lambda \unrhd \mu$.
\end{em}
\end{lemma}  

{\bf Proof:}

If the $\mu_{1}$ numbers of the first row of $t_{2}$ occur in different columns of $t_{1}$ then the latter has at least $\mu_{1}$ columns and so $\lambda_{1} \geq \mu_{1}$. Now consider the first two rows to deduce that $\lambda_{1} + \lambda_{2} \geq \mu_{1} + \mu_{2}$ and so on. $\Box$

The row-stabiliser $R_{t}$ of a tableau $t$ is the subgroup of $\Sigma_{n}$ which preserves each row of $t$ as a set - not element-wise! So for the left-hand tableau in the above example $R_{t}$  is 
\[  \Sigma( \{1,2,4,5 \} ) \times \Sigma( \{ 3,6,7 \} ) \times \Sigma( \{ 8 \} )    \]
where $\Sigma(S)$ denotes the group of permutations of the set $S$. 
\begin{definition}
\label{6.4}
\begin{em}

Two $\lambda$-tableaux $t_{1}$ and $t_{2}$ are equivalent, written $t_{1} \sim t_{2}$, if they belong to the same $\Sigma_{n}$-orbit (under the row-action). Each such orbit or equivalence class is called a tabloid, denoted by $\{ t \}$.

Now to define tabloid-orderings.

We say $\{ t_{1} \} < \{ t_{2} \}$ if and only if for some $i$ the following conditions hold:

(i) \  if $j > i$ then $j$ is in the same row of $\{ t_{1} \}$ and $\{ t_{2} \}$

(ii) \  $i$ is in a higher row of $\{ t_{1} \}$ than $ \{ t_{2} \}$.

Given a tableau $t$ let $m_{i,r}(t)$ denote the number of entries less than or equal to $i$ in the first $r$ rows of $t$. Then $\{ t_{1} \}  \unlhd \{ t_{2} \}$ if and only if for all $i$ and $r$ we have $m_{i,r}(t_{1}) \leq m_{i,r}(t_{2})$.

For tabloids $\{ t_{1} \} $ and $ \{ t_{2} \}$, we have that $\{ t_{1} \}  \lhd \{ t_{2} \}$ implies 
$\{ t_{1} \}  < \{ t_{2} \}$.

Suppose that $w < x$ such that $w$ is in the $a$-th row and $x$ in the $b$-th row of $t$. Then we characterise a tableau $t(w,x)$ by
\[   m_{i,r}(t(w,x))  - m_{i,r}(t) =  \left\{  \begin{array}{ll}
1 & {\rm if} \ b \leq r < a \ {\rm and} \ w \leq i < x, \\
\\
-1 & {\rm if} \ a \leq r < b \ {\rm and} \ w \leq i < x, \\
\\
0 & {\rm otherwise}.
\end{array} \right.      \]

This implies that $\{ t \}  \lhd  \{ t(w,x) \}$ if $w < x$ and $w$ is lower than $x$ in $t$.
\end{em}
\end{definition}
\begin{lemma}{(Adjacency in the $\lhd$ ordering  (\cite{GDJ} Lemma 3.16))}
\label{6.5}
\begin{em}

If $x-1$ is lower than $x$ in $t$, a $\lambda$-tableau, then there is no $\lambda$-tableau $t_{1}$ such that $ \{ t \} \lhd \{ t_{1} \}  \lhd  \{ t(x-1, x) \}$.
\end{em}
\end{lemma} 
\begin{definition}{(Specht modules)}
\label{6.6}
\begin{em}

Let $(\mu_{1}, \mu_{2}, \ldots, m_{k})$ be a partition of $n$. The row-stabiliser of the tabloid in which the rows are the sets of integers in the brackets of 
\[ (1, \ldots, \mu_{1})(\mu_{1}+1, \ldots , \mu_{1} + \mu_{2}) \ldots ( m_{1} + m_{2} + \ldots + m_{k-1} + 1, \ldots , n) \]
is the product subgroup
\[ \Sigma( \{ 1,  \ldots, \mu_{1} \}) \times  \Sigma( \{ \mu_{1}+1, \ldots , \mu_{1} + \mu_{2} \}) \times \ldots \]
Define $M^{\mu}$ to be the ${\mathbb Q}$-vector space whose basis consists of the 
$\frac{n!}{m_{1}! \mu_{2}! \ldots \mu_{k}!}$ tabloids in the $\Sigma_{n}$-orbit of the above $\mu$-tabloid. It is a cyclic ${\mathbb Q}[\Sigma_{n}]$-permutation module  generated by any tabloid equivalent to the above one.

In representation theory terms we have an isomorphism of ${\mathbb Q}[\Sigma_{n}]$-modules 
\[  M^{\mu}  \cong  {\rm Ind}_{\Sigma( \{ 1,  \ldots, \mu_{1} \}) \times   \ldots }^{\Sigma_{n}}( {\mathbb Q})  . \]

Suppose that $t$ is a tableau, defining the column stabiliser $C_{t}$ of $t$ in the manner analogous to that of $R_{t}$ we define 
\[   \kappa_{t} = \sum_{ \sigma \in C_{t}} \ {\rm sign}(g) \cdot g \in  {\mathbb Q}[\Sigma_{n}] . \]

The polytabloid associated to the tableau $t$ is given by
\[     e_{t} =  \kappa_{t}( \{ t \}) \in M^{\mu} . \]
The Specht module $S^{\mu}$ for the partition $\mu$ is the submodule of $M^{\mu}$ spanned by polytabloids.

Note that a polytabloid $e_{t}$ depends on the tableau $t$, not just on the tabloid $\{ t \}$.
All the tabloids involved in $e_{t}$ have coefficient equal to $\pm1$.
\end{em}
\end{definition}

Here is a small example. Suppose that 
\[ t =   \left(  \begin{array}{ccc}
2&5&1 \\
\\3&4& \\
\end{array} \right)    \]
and the columns of  $t$ are preserved by $1, (2,3), (4,5), (23)(45)$ so that 
\[ \kappa_{t} = 1 - (2,3) - (4,5) + (2,3)(4,5) \in {\mathbb Q}[\Sigma_{5}] . \]
If $\pi \in \Sigma_{5}$ then 
\[ \pi(t) =   \left(  \begin{array}{ccc}
\pi(2)&\pi(5)& \pi(1) \\
\\
\pi(3)& \pi(4)& \\
\end{array} \right)    \]
so that
 \[ \kappa_{\pi(t)} = 1 - (\pi(2),\pi(3)) - (\pi(4),\pi(5)) + (\pi(2),\pi(3))(\pi(4),\pi(5)) = \pi \kappa_{t} \pi^{-1} \in 
 {\mathbb Q}[\Sigma_{5}] . \]
 This property holds in general and proves that $S^{\mu}$ is a cyclic ${\mathbb Q}[\Sigma_{n}]$-module 
 generated by any polytabloid.
 
 If $\mu = (1^{n})$ then $S^{\mu} \cong {\mathbb Q}[\Sigma_{n}]$, the regular representation.
 
 \begin{lemma}{(\cite{GDJ} Lemma 4.6)}
 \label{6.7}
 \begin{em}
 
 Let $\lambda$ and $\mu$ be two partitions of $n$. Suppose that $t$ is a $\lambda$-tableau and $t'$ is a $\mu$-tableau and that $\kappa_{t}( \{ t' \}) = 0$. Then $\lambda \unrhd \mu$ and if $\lambda = \mu$ we have $\kappa_{t}( \{ t' \}) = \pm \kappa_{t}( \{ t \}) = \pm e_{t}$.
 
 Therefore, if $u \in M^{\mu}$ and $t$ is a $\mu$-tableau, then $\kappa_{t}(u)$ is a multiple of $e_{t}$.
 \end{em}
 \end{lemma} 
 \begin{definition}
 \label{6.8}
 \begin{em}
  Let $\langle - , - \rangle : M^{\mu} \times M^{\mu} \longrightarrow {\mathbb Q} $ be the symmetric
 ${\mathbb Q}[\Sigma_{n}]$-invariant bilinear form with respect to which the $\mu$-tabloids are an orthonormal basis.
 
 Importantly we have $\langle \kappa_{t}(u), v \rangle = \langle u, \kappa_{t}(v) \rangle$.
 \end{em}
 \end{definition} 
 \newpage
 
 \begin{theorem}{The Submodule Theorem (\cite{GDJ76})\footnote{A bit of old buffer's related anecdotage: In the 1980's I was the chair of the Canadian Mathematical Society's Research Committee. 
 The committee's job was to seek out and encourage research conferences in promising topics (established or potential) in order to influence (modernise?) mathematical research in Canada. As a result I learnt of and had to attend conferences on all sorts of areas of pure mathematics. Thus it was that one time I found myself at a meeting on $q$-series, listening to fascinating talks by the likes of Ian Macdonald, Jonathan Alperin, George Lusztig. At this meeting, in Edmonton, Alberta,  I learnt of this result and the following one by overhearing it mentioned in conversation by my friend Peter Hoffman. I immediately pricked up my ears since Gordon James$^{\dagger}$ had been another of my friends ever since the 1970's when we were both Cambridge dons. These are two elegant and important results!}}
 \label{6.9}
 \begin{em}
 
 If $U$ is a ${\mathbb Q}[\Sigma_{n}]$-submodule of $M^{\mu}$ then wither $S^{\mu} \subseteq U $ or $U  \subseteq (S^{\mu})^{\perp}$, the orthogonal complement of $S^{\mu}$.
 \end{em}
 \end{theorem} 
 
 {\bf Proof:}
 
 Suppose that $u \in U$ and $t$ is a $\mu$-tableau. Then $\kappa_{t}(u)$ is a multiple of $e_{t}$. If we can choose $u$ and $t$ such that the multiple is non-zero then $ e_{t} \in U$ and so $S^{\mu} \subseteq U$ since $e_{t}$ generates $S^{\mu}$. If, for every $t$ and $u$, $\kappa_{t}(u) = 0$ we have
 \[ 0 = \langle \kappa_{t}(u) , \{ t \} \rangle   = \langle u , \kappa_{t}(\{ t \}) \rangle =  \langle u , e_{t} \rangle     \]
 so that $U  \subseteq (S^{\mu})^{\perp}$ since $e_{t}$ generates $S^{\mu}$. $\Box$
 \begin{theorem}{(\cite{GDJ76} and (\cite{GDJ} \S11))}
 \label{6.10}
 \begin{em}
 
 $S^{\mu}/S^{\mu} \bigcap (S^{\mu})^{\perp}$ is either zero or absolutely irreducible \cite{JPS2004}. Furthermore, if it is non-zero, then $S^{\mu} \bigcap (S^{\mu})^{\perp}$ is the unique maximal submodule of $S^{\mu}$ and $S^{\mu}/S^{\mu} \bigcap (S^{\mu})^{\perp}$ is self-dual.
 
 In fact all the irreducible representations of $\Sigma_{n}$ are constructible in this manner.
 \end{em}
 \end{theorem}
 \begin{remark}
 \label{6.10c}
 \begin{em}
(i) \   Clearly $S^{(n)}/S^{(n)} \bigcap (S^{(n)})^{\perp} \cong {\mathbb Q}$, the trivial one-dimensional representation and $S^{(1^{n})}/S^{(1^{n})} \bigcap (S^{(1^{n})})^{\perp} \cong {\mathbb Q}_{{\rm sign}}$, the one-dimensional representation on which $\Sigma_{n}$ acts via the sign character and $S^{(1^{n})} = {\mathbb Q}[\Sigma_{n}]$.

(ii)  \  For ${\mathbb Q}$-representations of $\Sigma_{n}$ we have $S^{\mu} \bigcap (S^{(\mu})^{\perp} = 0$ so that the ${\mathbb Q}$-Specht module $S^{\mu}$ gives an absolutely irreducible ${\mathbb Q}[ \Sigma_{n}]$-module. We shall denote the irreducible complex representation of $\Sigma_{n}$ given on $S^{\mu} \otimes_{{\mathbb Q}} {\mathbb C}$ by $\rho_{\mu}$ where $\mu$ is a partition of $n$.
 
 (iii) If $\lambda'$ is the conjugate of $\lambda$ then, as ${\mathbb Q}$-representations, 
 \[   S^{(\lambda')}   \cong S^{(\lambda)} \otimes {\mathbb Q}_{{\rm sign}} , \]
 as one would naturally guess.  The result 
(\cite{GDJ} Theorem 6.7)  gives a non-zero homomorphism between these two ${\mathbb Q}[\Sigma_{n}]$-modules, which must be an isomorphism by irreducibility. Alternatively this isomorphism follows directly from the definition of the ${\mathbb Q}$-Specht module.
 \end{em}
 \end{remark} 
  \begin{definition}
 \label{6.11}
 \begin{em}
 
 $t$ is a standard tableau if the numbers increase along the rows and down the columns of $t$. Then $\{ t \}$ is called a standard tabloid if there is a standard tableau in the equivalence class (i.e. orbit) 
 $\{ t \}$ in which case $e_{t}$ is called a standard polytabloid.
 
 A standard tabloid contains a unique standard tableau since the rows are increasing and the permutation action only permutes each row. However a standard polytabloid may involve more than one standard tabloid (see (\cite{GDJ} Example 5.2).
 \end{em}
 \end{definition}
 \begin{theorem}{Basis Theorem for $S^{\mu}$ (\cite{GDJ} Theorem 8.4)}
 \label{6.12}
 \begin{em}
 
 The set $ \{ e_{t} \ | \ t \ {\rm is \ a \ standard} \  \mu-{\rm tableau} \}$ is a ${\mathbb Q}$-basis for $S^{\mu}$
 \end{em}
 \end{theorem}
 
 \begin{theorem}{(The Branching Theorem)} (\cite{GDJ} Theorem 9.2 p.34)
 \label{6.13}
 \begin{em}
 
 (i)   \ ${\rm Ind}_{\Sigma_{n}}^{\Sigma_{n+1}}( S^{\mu}) \cong \oplus_{\lambda} \ S^{\lambda} $, the direct sum over all the diagrams $[\lambda]$ obtained by adding a node to the diagram $[\mu]$.
 
 (ii) \  ${\rm Res}_{\Sigma_{n-1}}^{\Sigma_{n}}( S^{\mu})  \cong \oplus_{\lambda} \ S^{\lambda} $, the direct sum over all the diagrams $[\lambda]$ obtained by deleting a node to the diagram $[\mu]$.
 \end{em}
 \end{theorem}
 
  For example in (i), if $[\mu] = [(4, 2^{2}),1)]$ then the $[\lambda]$'s are
  $[(5, 2^{2}, 1)]$, $[(4,3,2,1)]$, $[(4, 2^{3})]$ and $[(4, 2^{2}, 1^{2})]$ whereas,
 for example in (ii), if $[\mu] = [(4, 2^{2}),1)]$ then the $[\lambda]$'s are
 $ [(3, 2^{2}),1)]$, $ [(4, 2, 1^{2}))]$ and $ [(4, 2^{2}))]$.
\vspace{10pt} 

 Before I begin the study of numerical invariants let us pause for an example. Here is the $\Sigma_{5}$ character table (\cite{GDJ} p.24) - rows are the Specht module 
 ${\mathbb Q}$-representations and columns are their character valued at conjugacy classes given by the partitions of cycle decompositions.
 
 \[ \begin{array}{|c|ccccccc|}
 \hline
 \hline
  & (5) & (4,1) & (3,2) & (3,1^{2}) & (2^{2}, 1) & (2, 1^{3}) & (1^{5}) \\
  \hline
 & &&&&&& \\
  (5) & 1&1&1&1&1&1&1 \\
& &&&&&&  \\
  (4,1) & -1 &0&-1&1&0&2&4 \\
& &&&&&&  \\
  (3,2) & 0&-1&1&-1&1&1&5 \\
& &&&&&&  \\
  (3,1^{2}) & 1&0&0&0&-2&0&6\\
& &&&&&&  \\
(2^{2}, 1) & 0&1&-1&-1&1&-1&5 \\
& &&&&&&\\
(2, 1^{3}) & -1&0&1&1&0&-2&4 \\
& &&&&&&\\
(1^{5}) &  1&-1&-1&1&1&-1&1 \\
  \hline
 \hline
 \end{array} \]
 \begin{definition}{Symmetric group analogues of Kondo-Gauss sums}
\label{6.14}
\begin{em}

Let $\lambda$ be a partition of $n$. The length of $\lambda$ is simply the number of strictly positive integers in the partition. If $\sigma  \in \Sigma_{n}$, the symmetric group, then the cycle decomposition of $\sigma$ gives a partition of $n$ and we define the length of $\sigma$, denoted by $l(\sigma)$, to be the length of that partition.

Let $x$ be a complex number. For $\sigma \in \Sigma_{n}$ define a measure $\psi_{x}$\footnote{Warning: In the section on the general wreath product we shall have to keep track of which symmetric group $\Sigma_{n}$ we are concerned with and at that point $\psi_{x}$ is apt to become $\psi_{n,x}$.} on $\Sigma_{n}$ by
\[ \psi_{x}(\sigma) =   x^{ l(\sigma)} .  \]
The measure $\psi_{x}$ can be used to give integrals of characters of representations of the symmetric group, in the spirit of Kondo-Gauss sums. This function is not an additive character - but it is additive on pairs of disjoint permutations, making it very suitable for the PSH algebra of $\Sigma_{n} $.

If $H \subseteq \Sigma_{n}$ and $\rho$ is a complex representation of $H$ define
\[   W_{H}^{x}(\rho) = \frac{1}{{\rm dim}_{{\mathbb C}}(\rho) } \ \sum_{h \in H} \ \chi_{\rho}(h) \psi_{x}(h) .  \]
\end{em}
\end{definition}

In the special case when $H = \Sigma_{n}$ and $\rho = \rho_{\lambda}$ in the notation of Remark \ref{6.10c}(ii) we have the following formula:
\begin{theorem}{(F. Mezzadri 2010; appearing in \cite{KMS2011}; see also (\cite{Sn18} pp.307-8))}
\label{6.15}
\begin{em}
\[    W_{\Sigma_{n}}^{x}( \rho_{\lambda}) =  \prod_{i=1}^{{\rm length}(\underline{\lambda}) }  \  
(x - i+1)(x-i+2) \ldots (x + \lambda_{i} - i) =  \prod_{(i,j) \in [ \lambda ]}  (x - i + j).\] 

In \cite{KMS2011} this monic polynomial with integer coefficients is denoted by $ f_{\underline{\lambda}}(x)$.
\end{em}
\end{theorem} 

This is the integral to which I referred in the historical footnote at the beginning of this section. When one is familiar with the material in \cite{GDJ}, which I have just recapitulated in this section, it is an easy induction. The induction and the behaviour of  $W_{\Sigma_{n}}^{x}(\rho) $ with respect to the PSH algebra product (\cite{Sn18} Appendix III Proposition 4.3) will occupy the rest of this section\footnote{The next section shows how all these formulae generalise to the case of the PSH algebra of wreath products.}.
\begin{lemma}
\label{6.16}
\begin{em}

(i)  \ If $\lambda'$ is the conjugate of $\lambda$ then
\[   W_{\Sigma_{n}}^{x}( \rho_{\lambda}) = (-1)^{n} W_{\Sigma_{n}}^{-x}( \rho_{\lambda'}) .\]

(ii) \ Also we have
\[ \prod_{(i,j) \in [ \lambda ]}  (x - i + j) = (-1)^{n}  \prod_{(i,j) \in [ \lambda' ]}  (- x - i + j) . \]
\end{em}
\end{lemma} 

{\bf Proof}

Part (i) follows from the fact that $\chi_{\rho_{\lambda'} }(\sigma) = \chi_{\rho_{\lambda'} }(\sigma) \cdot {\rm sign}(\sigma) $ and the sign of $\sigma$ is minus one to the number of even length cycles in the cycle decomposition of $\sigma$. Part (ii) is clear from the definition of the conjugate diagram in Definition \ref{6.2}. $\Box$
\begin{dummy}
\label{6.18}
\begin{em}

Let $V_{1}, V_{2}$ be complex representations of $\Sigma_{k}, \Sigma_{n-k}$ respectively and let $m(V_{1} \otimes V_{2}$ be their product in the PSH algebra of \S4. Recall that this product is the
induced representation up to $\Sigma_{n}$.

 The trace formula for the character of this representation of $\Sigma_{n}$ is given by 
\[  \chi_{m(V_{1} \otimes V_{2})}(\hat{\sigma})  = \frac{1}{k!(n-k)!} \sum_{\tau \in \Sigma_{n}, \tau \hat{\sigma} \tau^{-1} \in \Sigma_{k} \times \Sigma_{n-k}} \  \chi_{ V_{1} \otimes V_{2}}( \tau \hat{\sigma} \tau^{-1} )   \]
which is zero unless $\hat{\sigma}$ is conjugate in $\Sigma_{n}$ to an element $(\sigma, \sigma') \in 
\Sigma_{k} \times  \Sigma_{n-k}$. 

Two elements of $\Sigma_{k} \times \Sigma_{n-k}$ are conjugate in $\Sigma_{n}$ if and only if they are conjugate in $\Sigma_{k} \times \Sigma_{n-k}$, because cycle shape determines conjugacy.

Suppose that $\hat{\sigma}$ is conjugate in $\Sigma_{n}$ to an element $(\sigma, \sigma') \in 
\Sigma_{k} \times  \Sigma_{n-k}$.  Let $r_{i}, r'_{i}$ be the numbers of $i$-cycles in the cycle decomposition for $\sigma, \sigma'$ respectively. The number of distinct elements in the $\Sigma_{n}$-conjugacy class of $(\sigma, \sigma')$ is
\[   \frac{n!}{1^{r_{1} + r'_{1}}(r_{1} + r'_{1})! 2^{r_{2} + r'_{2}}(r_{2} + r'_{2})! \ldots  }    \]
while the number in the $\Sigma_{k}$-conjugacy class of $\sigma$ is
\[   \frac{k!}{1^{r_{1} }(r_{1} )! 2^{r_{2} }(r_{2})! \ldots  }    \]
and  the number in the $\Sigma_{n-k}$-conjugacy class of $\sigma'$ is
\[   \frac{(n-k)!}{1^{r'_{1}}(r'_{1})! 2^{ r'_{2}}( r'_{2})! \ldots  } .   \]
Therefore the denominators are the orders of the relevant centralisers.

Next we want to simplify the formula for $ \chi_{m(V_{1} \otimes V_{2})}(\hat{\sigma}) $. If $\tau \hat{\sigma} \tau^{-1} \in \Sigma_{k} \times \Sigma_{n-k}$ there exists $\mu \in \Sigma_{k} \times \Sigma_{n-k}$ such that $\mu \tau \hat{\sigma} \tau^{-1} \mu^{-1} = \hat{\sigma}$. Setting $\tau_{1} = \mu \tau$
we see that $\tau_{1} \in Z_{\Sigma_{n}}(\hat{\sigma})$, the centraliser of $\hat{\sigma}$ in $\Sigma_{n}$.
Therefore $\tau = \mu^{-1} \tau_{1}$ and we have a surjective map
\[  (\Sigma_{k} \times \Sigma_{n-k}) \times Z_{\Sigma_{n}}(\hat{\sigma})  \longrightarrow 
\{ \tau \in \Sigma_{n}, \tau \hat{\sigma} \tau^{-1} \in \Sigma_{k} \times \Sigma_{n-k} \} \]
sending $(\nu , \tau_{1}) $ to $\nu \tau_{1}$. Also $\nu \tau_{1} = \nu' \tau'_{1}$ if and only if
\[ ( \nu')^{-1} \nu = \tau'_{1} \tau_{1}^{-1} \in Z_{\Sigma_{k} \times \Sigma_{n-k}}(\sigma, \sigma').\]
Therefore if $\lambda = ( \nu')^{-1} \nu$ then $(\nu, \tau_{1}) = (\nu' \lambda, \lambda^{-1} \tau')$. Therefore there is a bijection
\[  (\Sigma_{k} \times \Sigma_{n-k}) \times_{Z_{\Sigma_{k} \times \Sigma_{n-k}}(\sigma, \sigma')} Z_{\Sigma_{n}}(\hat{\sigma})  \leftrightarrow 
\{ \tau \in \Sigma_{n}, \tau \hat{\sigma} \tau^{-1} \in \Sigma_{k} \times \Sigma_{n-k} \}.  \]
Therefore the number of $\tau$'s in the formula for $ \chi_{m(V_{1} \otimes V_{2})}(\hat{\sigma}) $ is
\[ \frac{k!(n-k)! 1^{r_{1} + r'_{1}}(r_{1} + r'_{1})! 2^{r_{2} + r'_{2}}(r_{2} + r'_{2})! \ldots }{1^{r_{1} }(r_{1} )! 2^{r_{2} }(r_{2})! \ldots  1^{r'_{1}}(r'_{1})! 2^{ r'_{2}}( r'_{2})! \ldots  }  \]
and each $\tau \hat{\sigma} \tau^{-1}$ gives an element in the $\Sigma_{k} \times \Sigma_{n-k}$-conjugacy class of $\hat{\sigma} = (\sigma, \sigma')$ we find that
\[  \chi_{m(V_{1} \otimes V_{2})}(\hat{\sigma})  =   \frac{ 1^{r_{1} + r'_{1}}(r_{1} + r'_{1})! 2^{r_{2} + r'_{2}}(r_{2} + r'_{2})! \ldots }{1^{r_{1} }(r_{1} )! 2^{r_{2} }(r_{2})! \ldots  1^{r'_{1}}(r'_{1})! 2^{ r'_{2}}( r'_{2})! \ldots  }   \chi_{ V_{1}}(\sigma)  \chi_{V_{2}}(\sigma' )   \]
if $\hat{\sigma}$ is $\Sigma_{n}$-conjugate to $(\sigma, \sigma') \in \Sigma_{k} \times \Sigma_{n-k}$ and zero otherwise.
\end{em}
\end{dummy}
\begin{theorem}{(\cite{Sn18} Appendix III Proposition 4.3)}
\label{6.19}
\begin{em}

If $V_{1}$ and $V_{2}$ are representations of $\Sigma_{k}$ and $\Sigma_{n-k}$ respectively and 
\linebreak
$m(V_{1} \otimes V_{2})$ is the product in the PSH algebra of \S4 then
\[     W_{\Sigma_{n}}^{x}( m(V_{1} \otimes V_{2})) =    W_{\Sigma_{k}}^{x}( V_{1})  \cdot 
W_{\Sigma_{n-k}}^{x}( V_{2})  .\]
\end{em}
\end{theorem}

{\bf Proof:} 

 From the preceding discussion about the character values of 
\[ {\rm Ind}_{\Sigma_{k} \times  \Sigma_{n-k}}^{\Sigma_{n}}(V_{1} \otimes V_{2}) = m(V_{1} \otimes V_{2}) \]
we see that 
\[ \begin{array}{l}
   W_{\Sigma_{n}}^{x}( m(V_{1} \otimes V_{2})) \\
  \\
  =   \frac{k!(n-k)!}{n! {\rm dim}_{ {\mathbb C}}(V_{1}){\rm dim}_{{\mathbb C}}(V_{2}) } \sum_{\hat{\sigma} \in \Sigma_{n}} \    \chi_{m(V_{1} \otimes V_{2})}(\hat{\sigma}) \psi_{x}(\hat{\sigma}) \\
  \\
  =   \frac{k!(n-k)!}{n! {\rm dim}_{{\mathbb C}}( V_{1}){\rm dim}_{{\mathbb C}}(V_{2}) } \sum_{\hat{\sigma} \in \Sigma_{n}} \     \frac{ 1^{r_{1} + r'_{1}}(r_{1} + r'_{1})! 2^{r_{2} + r'_{2}}(r_{2} + r'_{2})! \ldots }{1^{r_{1} }(r_{1} )! 2^{r_{2} }(r_{2})! \ldots  1^{r'_{1}}(r'_{1})! 2^{ r'_{2}}( r'_{2})! \ldots  }   \chi_{ V_{1}}(\sigma) \\
  \\
  \hspace{160pt} \times   \chi_{V_{2}}(\sigma' )  \psi_{x}(\hat{\sigma}) 
\end{array} \]
where the sum is taken over those $\hat{\sigma}$ in $\Sigma_{n}$ which are conjugate to some $(\sigma, \sigma')$ in $\Sigma_{k} \times \Sigma_{n-k}$. Now if $\hat{\sigma} = (\sigma, \sigma')$ with lengths as in the preceding discussion then the ratio of conjugacy class sizes satisfies
\[ \begin{array}{l}
\frac{|\{ \Sigma_{n}-{\rm conjugates  \ of}  \ (\sigma, \sigma')\}|}{|\{ \Sigma_{k \times \Sigma_{n-k}}-{\rm conjugates  \ of } \ (\sigma, \sigma')\}|} \\
\\
=    \frac{n!}{1^{r_{1} + r'_{1}}(r_{1} + r'_{1})! 2^{r_{2} + r'_{2}}(r_{2} + r'_{2})! \ldots  }  
   \frac{1^{r_{1} }(r_{1} )! 2^{r_{2} }(r_{2})! \ldots  }{k!}   
  \frac{1^{r'_{1}}(r'_{1})! 2^{ r'_{2}}( r'_{2})! \ldots  }{(n-k)!}    
\end{array} \]
Therefore we may re-write the function as a sum over $\Sigma_{k} \times \Sigma_{n-k}$ in the form

\[  \begin{array}{l}
   W_{\Sigma_{n}}^{x}( m(V_{1} \otimes V_{2}))   \\
 \\
 =    \frac{1}{ {\rm dim}_{{\mathbb C}}( V_{1}) {\rm dim}_{{\mathbb C}}(V_{2}) }  \sum_{(\sigma, \sigma') \in \Sigma_{k} \times \Sigma_{n-k} }  \chi_{ V_{1}}(\sigma)  \chi_{V_{2}}(\sigma' )   \psi_{x}(\hat{\sigma}) \\
 \\
 =  W_{\Sigma_{k}}^{x}( V_{1})  \cdot  W_{\Sigma_{n-k}}^{x}( V_{2})  ,
\end{array}  \]
as required, since $ \psi_{x}(\hat{\sigma}) =  \psi_{x}(\sigma) \psi_{x}(\sigma')$. $\Box$
\begin{corollary}{$_{}$}
\label{6.20}
\begin{em}

If ${\bf 1}$ denotes the trivial one-dimensional representation then 
\linebreak
${\rm Ind}_{\Sigma_{n}}^{\Sigma_{n+1}}(V) = m(V \otimes {\bf 1})$ so that $W_{\Sigma_{n+1}}^{x}({\rm Ind}_{\Sigma_{n}}^{\Sigma_{n+1}}(V) ) = x  W_{\Sigma_{n}}^{x}(V )$
\end{em}
\end{corollary} 
\begin{dummy}{Proof of Theorem \ref{6.15}}
\label{6.20}
\begin{em}

The proof is by induction, using the Branching Rule of Theorem \ref{6.13}(i). Recall that $\rho_{\mu}$ denotes the representation of $\Sigma_{n}$ given by the complexification of the Specht module $S^{\mu}$.
\end{em}
\end{dummy} 

The branching rule and the analogue of Lemma \ref{5.3}  yield the relation
\[ \begin{array}{l}
(n+1) {\rm dim}(S^{\mu}) x \cdot  W_{\Sigma_{n}}^{x}(\rho_{\mu})    \\
\\
= {\rm dim}(  {\rm Ind}_{\Sigma_{n}}^{\Sigma_{n+1}}( \rho_{\mu} ) ) 
W_{\Sigma_{n+1}}^{x}({\rm Ind}_{\Sigma_{n}}^{\Sigma_{n+1}}( \rho_{\mu}) )  \\
\\
= \sum_{\lambda}  \ {\rm dim}(S^{\lambda})   W^{x}_{\Sigma_{n+1}}(\rho_{\lambda})  .
\end{array} \] 
where the sum is over all the diagrams $[\lambda]$ obtained by adding a node to the diagram $[\mu]$.

Following (\cite{KMS2011} and \cite{Sn18} Appendix III \S4) denote the polynomial on the right-hand side of the equation of Theorem \ref{6.15} by $ f_{ S^{\mu} }(x) $ for the representation $\rho_{\mu}$.

Furthermore, when $\lambda$ is obtained by adding the $(i(\lambda),j(\lambda))$-node to the diagram $[\mu]$ to obtain $[\lambda]$, we have the relation
\[  f_{ S^{\lambda} }(x) = f_{ S^{\mu} }(x) \cdot (x - i(\lambda) + j(\lambda)) . \]
Therefore, if we replace each $ W^{x}_{\Sigma_{n+1}}(\rho_{\lambda}) $ in the sum of $\lambda$'s by $ f_{ S^{\lambda} }(x) $ we obtain
\[  \begin{array}{l}
\sum_{\lambda}  \ {\rm dim}(S^{\lambda})  f_{ S^{\mu} }(x) (x - i(\lambda) + j(\lambda)) \\
\\
=   (n+1) {\rm dim}(S^{\mu}) x \cdot   f_{ S^{\mu} }(x)  -   f_{ S^{\mu} }(x) \cdot ( \sum_{\lambda}  \ {\rm dim}(S^{\lambda})  (i(\lambda) - j(\lambda))) .
\end{array}   \]
The sum
\[  \sum_{\lambda}  \ {\rm dim}(S^{\lambda})  (i(\lambda) - j(\lambda))  \]
 is easily seen to be zero. The dimensions in the sum are unchanged by applying the conjugation involution $\mu \mapsto \mu'$. However conjugation interchanges the $i$'s and $j$'s, as the $(i,j)$'s are precisely those where an extra node can be added to the diagram for $S^{\mu}$ and the $(j,i)$'s are precisely those where an extra node can be added to the diagram for $S^{\mu'}$,  so the expression must be zero. 
 
 Therefore, if Theorem \ref{6.15} holds for all the $\lambda$ in the branching rule, then it holds for $\mu$. 
 This is the basis for an induction which starts with the trivial cases $(n)$ and $(1^{n})$ for all $n$.
 $\Box$
 
 The proof of Theorem \ref{5.2} also yields the following result.
 \begin{theorem}{(Analogue of Theorem \ref{5.2})}
\label{6.21}
\begin{em}

Let $\sigma $ be a finite-dimensional representation of $H \subseteq \Sigma_{n}$. Then for any subgroup $J$ such that $H \subseteq J  \subseteq \Sigma_{n}$
\[   W_{H}^{x}(\sigma) = W_{J}^{x}( {\rm Ind}_{H}^{J}(\sigma)) .\]
\end{em}
\end{theorem}
 
 \section{The general wreath product}

In this section we consider the wreath-product $\Sigma_{n} \int H$ which is a subgroup of $GL_{nd}{\mathbb F}_{q}$ if $H$ is subgroup of $GL_{d}{\mathbb F}_{q}$. We start by extending, in a preliminary way, the definition of the symmetric group numerical invariant of \S6.

\begin{definition}
\label{7.1}
\begin{em}
 Suppose that $H \subseteq GL_{m}{\mathbb F}_{q}$ ($q = p^{d}$, $p$ prime) and that $\rho$ is a finite-dimensional complex representation of the subgroup $J$ of the wreath product $\Sigma_{n} \int H$ define, for $x \in {\mathbb C}$,
\[ W_{n, J}^{x, H}(\rho) = \frac{1}{ {\rm dim}_{{\mathbb C}}(\rho) }   
 \sum_{ X = (\sigma , \alpha_{1}, \ldots , \alpha_{n})  \in J \subseteq  \Sigma_{n} \int H} 
   \chi_{\rho}(X)  \psi_{n,x}(\sigma) \Lambda_{n, x, H}(  (\sigma , \alpha_{1}, \ldots , \alpha_{n})) \]
 where
   \[ \Lambda_{n,x, H}(  (\sigma , \alpha_{1}, \ldots , \alpha_{n}))  = \prod_{\stackrel{  {\rm cycles} \ ( x_{1}, x_{2}, \ldots , x_{r}) }{ \rm  in \   \sigma \in \Sigma_{n}} }    \    e^{  2 \pi \sqrt{-1} {\rm Trace}_{ {\mathbb F}_{q}/ {\mathbb F}_{p}} ({\rm Trace}_{GL_{m}}(\alpha_{x_{1}} \alpha_{x_{2}} \ldots \alpha_{x_{r}}) )/p    } .      \]
\end{em}
\end{definition} 

\begin{theorem}{$_{}$}
\label{7.2}
\begin{em}
 If $G \subset J \subseteq \Sigma_{n} \int H$ in the situation of Definition \ref{7.1} and $\lambda$ is a complex representation of $G$ then
    \[ \begin{array}{l}
 W_{n, J}^{x, H}({\rm Ind}_{G}^{J}(\lambda) ) \\
\\
 =   \frac{1}{|J| {\rm dim}_{{\mathbb C}}(\lambda) }   
 \sum_{  \stackrel{W = (\sigma , \alpha_{1}, \ldots , \alpha_{n})  \in G  }{ Y = (\sigma' , \alpha'_{1}, \ldots , \alpha'_{n})  \in J   } \ Y^{-1}WY \in J  } 
   \chi_{\lambda}(W)  \psi_{n,x}(\sigma) \mu_{n, x, H}(W,Y) 
   \end{array}   \]
   where
   \[ \mu_{n, x, H}(W,Y) =    \prod_{\stackrel{  {\rm cycles} \ ( x_{1}, x_{2}, \ldots , x_{r}) }{ \rm  in \   (\sigma')^{-1} \sigma \sigma'}  }    \    e^{  2 \pi \sqrt{-1} {\rm Trace}_{ {\mathbb F}_{q}/ {\mathbb F}_{p}} ({\rm Trace}_{GL_{m}}(\alpha_{x_{1}} \alpha_{x_{2}} \ldots \alpha_{x_{r}}) )/p    } \]
\end{em}
\end{theorem}

{\bf Proof:}

If $\rho = {\rm Ind}_{G}^{J}(\lambda)$ then
\[ \begin{array}{l}
W_{n, J}^{x, H}(\rho) \\
\\
= \frac{|G|}{|J| {\rm dim}_{{\mathbb C}}(\lambda) }   
 \sum_{ X = (\sigma , \alpha_{1}, \ldots , \alpha_{n})  \in J \subseteq  \Sigma_{n} \int H} 
   \chi_{\rho}(X)  \psi_{n,x}(\sigma) \Lambda_{n, x, H}( X) \\
   \\
 =   \frac{1}{|J| {\rm dim}_{{\mathbb C}}(\lambda) }   
 \sum_{  \stackrel{X = (\sigma , \alpha_{1}, \ldots , \alpha_{n})  \in J  }{ Y = (\sigma' , \alpha'_{1}, \ldots , \alpha'_{n})  \in J   } \ YXY^{-1} \in G  } 
   \chi_{\lambda}(YXY^{-1})  \psi_{n,x}(\sigma) \Lambda_{n, x, H}(  X) .
   \end{array}   \]
There is a bijection of sets of the form
\[  \begin{array}{c}
  \{  (X, Y) \in J \times J \ | \ YXY^{-1} \in G \}    \\
  \\
  \updownarrow   \\
  \\
    \{  (W, Y) \in G \times J  \}  .
  \end{array}  \]
  The downwards map is given by $(X,Y) \mapsto (YXY^{-1}, Y)$ and the upwards map is given by
  $(W,Y) \mapsto  (Y^{-1}WY, Y)$. Therefore we have 
  \[ \begin{array}{l}
W_{n, J}^{x, H}(\rho) \\
\\
 =   \frac{1}{|J| {\rm dim}_{{\mathbb C}}(\lambda) }   
 \sum_{  \stackrel{W = (\sigma , \alpha_{1}, \ldots , \alpha_{n})  \in G  }{ Y = (\sigma' , \alpha'_{1}, \ldots , \alpha'_{n})  \in J   } \ Y^{-1}WY \in J  } 
   \chi_{\lambda}(W)  \psi_{n,x}(\sigma) \Lambda_{n, x, H}(  Y^{-1}WY) 
   \end{array}   \]
   because $ \psi_{n,x}((\sigma')^{-1} \sigma \sigma') =  \psi_{n,x}(\sigma)$.
  
  Recall that the product in the wreath product is given by
   \[ (\sigma , \alpha_{1}, \ldots , \alpha_{n}) \cdot (\sigma' , \alpha'_{1}, \ldots , \alpha'_{n}) =
( \sigma \sigma' , \alpha_{1} \alpha'_{\sigma(1)} , \alpha_{2} \alpha'_{\sigma(2)} , \ldots , \alpha_{n} \alpha'_{\sigma(n)}) \]
so that 
\[ (\sigma'  , \alpha'_{1}, \ldots   , \alpha'_{n})^{-1} = ( (\sigma')^{-1}  , (\alpha'_{(\sigma')^{-1}(1)})^{-1}, \ldots , ( \alpha'_{(\sigma')^{-1}(n)})^{-1})  . \]
Therefore we have
\[ \begin{array}{l}
Y^{-1}WY \\
\\
= ( (\sigma')^{-1}  , (\alpha'_{(\sigma')^{-1}(1)})^{-1}, \ldots ) \dot ((\sigma , \alpha_{1}, \ldots ) \cdot (\sigma' , \alpha'_{1}, \ldots ))  \\
\\
= ( (\sigma')^{-1}  , (\alpha'_{(\sigma')^{-1}(1)})^{-1}, \ldots ) \dot  
( \sigma \sigma' , \alpha_{1} \alpha'_{\sigma(1)} , \alpha_{2} \alpha'_{\sigma(2)} , \ldots ) \\
\\
= (( \sigma')^{-1}  \sigma \sigma' ,  (\alpha'_{(\sigma')^{-1}(1)})^{-1}  \alpha_{ (\sigma')^{-1}(1)} \alpha'_{  (\sigma')^{-1} \sigma(1)} , \ldots )  .
\end{array} \]
Suppose that $(j_{1}, j_{2}, \ldots , j_{t})$ is a cycle in the cycle decomposition of $\sigma$. If $(\sigma')^{-1}(j_{s}) = k_{s}$ then $(\sigma')^{-1}  \sigma \sigma'(k_{s}) = (\sigma')^{-1}(j_{s+1}) = k_{s+1}$ (with the convention that $k_{t+1} = k_{1}$) so that $(k_{1}, \ldots , k_{t})$ is a cycle in 
$ (\sigma')^{-1}  \sigma \sigma'$ and we have
\[  (\alpha'_{(\sigma')^{-1}(j_{s})})^{-1}  \alpha_{ (\sigma')^{-1}(j_{s})} \alpha'_{  (\sigma')^{-1} \sigma(j_{s})}=   (\alpha'_{k_{s}})^{-1}  \alpha_{ k_{s}} \alpha'_{ k_{s+1}} .\]
Therefore the product of these terms, index by the cycle of $k_{i}$'s in cyclic order is equal to
\[     (\alpha'_{k_{1}})^{-1}  \alpha_{ k_{1}} \alpha_{k_{2}} \ldots  \alpha_{k_{t}}\alpha'_{ k_{1}}    \]
and
\[ \begin{array}{ll}
{\rm Trace}_{GL_{m}}( (\alpha'_{k_{1}})^{-1}  \alpha_{ k_{1}} \alpha_{k_{2}} \ldots  \alpha_{k_{t}} \alpha'_{ k_{1}} ) & =  {\rm Trace}_{GL_{m}}(  \alpha_{ k_{1}} \alpha_{k_{2}} \ldots  \alpha_{k_{t}}) \\
\\
& =   {\rm Trace}_{GL_{m}}(  \alpha_{ (\sigma')^{-1}(j_{1})} \ldots  \alpha_{(\sigma')^{-1}(j_{t})}) . 
\end{array} \]
This observation completes the proof. $\Box$

\begin{theorem}{$_{}$}
\label{7.3}
\begin{em}

If $H$ is a subgroup of $GL_{m}{\mathbb F}_{q}$  and $S^{\lambda}$ is the finite-dimensional complex representation of $\Sigma_{n}$ given by the Specht module of the partition $\lambda$ then
 \[  
 W^{x,H} _{n,\Sigma_{n} \int H}( {\rm Ind}_{\Sigma_{n}}^{\Sigma_{n} \int H}(S^{\lambda})) =    f_{\lambda}(x e^{  \frac{2 \pi \sqrt{-1} md }{p}    })  . \]
 
 Here $f_{\lambda}(x)$ is the polynomial associated to the $\Sigma_{n}$-representation $S^{\lambda}$ which appeared in \S\ref{6.20} (the proof of Theorem \ref{6.15}).
\end{em}
\end{theorem}

{\bf Proof:}

 In the formula of Theorem \ref{7.2} every $\alpha_{i} = 1$. Therefore for each cycle, of length $t$ say, we contribute a factor $x e^{  2 \pi \sqrt{-1} md/p    }  $ to $ \psi_{n,x}(\sigma) \Lambda_{n, x, H}(  Y^{-1}WY) $ which multiply together to give $ \psi_{n,x e^{  2 \pi \sqrt{-1} md/p    } }(\sigma)$. We receive the same contribution for each $Y$ which cancels out the $|J|$ in the denominator. $\Box$ 
 
 Theorem \ref{7.2} shows that $ W_{n, J}^{x, H}(\rho ) $ is not as inductive as the classical Kondo-Gauss sum because the measure is not sufficiently invariant to ensure that  $W_{n, J}^{x, H}({\rm Ind}_{G}^{J}(\lambda) ) $ is equal to $ W_{n, G}^{x, H}(\lambda )$.
 
 To emphasise this problem we shall pause for a small explicit example, which also shows that $ W_{n, J}^{x, H}(\rho ) $ does not respect the wreath product PSH algebra structure.

\begin{example}
\label{7.4}
\begin{em}

 $G = \Sigma_{2} \propto \{ (\alpha_{1}, \alpha_{2}, 1)\}$ and $J = \Sigma_{3} \int H$.
Therefore we have $X = (1, \alpha_{1}, \alpha_{2}, 1))$ and $X = ((1,2), \alpha_{1}, \alpha_{2}, 1)$ which $Y = (\sigma', \alpha'_{1}, \alpha'_{2}, \alpha'_{3})$ with $\sigma'$ in $\{1, (1,2), (1,3), (2,3), 123), 132) \}$.

When $\sigma = 1$ we have
\[ \begin{array}{|c|c|c|}
\hline
\hline
\sigma & \sigma' & Y^{-1}XY \\
\hline
\hline
1 & 1 &  (( 1,  (\alpha'_{1})^{-1}  \alpha_{ (1} \alpha'_{ 1} , (\alpha'_{2})^{-1}  \alpha_{ 2} \alpha'_{ 2} ,(\alpha'_{3})^{-1}  \alpha_{3} \alpha'_{ 3} )   \\
\hline
1 &  (1,2) &  ((1 ,  (\alpha'_{2})^{-1}  \alpha_{ 2} \alpha'_{ 2} , (\alpha'_{1})^{-1}  \alpha_{ 1} \alpha'_{ 1} , ( \alpha'_{3})^{-1}  \alpha_{ 3} \alpha'_{ 3} )   \\
\hline
1 & (1,3) &  ((1 ,  (\alpha'_{3})^{-1}  \alpha_{ 3} \alpha'_{ 3} , (\alpha'_{2})^{-1}  \alpha_{ 2} \alpha'_{ 2} , (\alpha'_{1})^{-1}  \alpha_{ 1} \alpha'_{ 1} )   \\
\hline
1 &  (2,3) &  ((1 ,  (\alpha'_{1})^{-1}  \alpha_{ 1} \alpha'_{ 1} ,( \alpha'_{3})^{-1}  \alpha_{ 3} \alpha'_{ 3} ,(\alpha'_{2})^{-1}  \alpha_{ 2} \alpha'_{ 2} )   \\
\hline
1 &  (1,2,3) &  ((1,  (\alpha'_{3})^{-1}  \alpha_{ 3} \alpha'_{ 3} , (\alpha'_{1})^{-1}  \alpha_{ 1} \alpha'_{  1} , (\alpha'_{2})^{-1}  \alpha_{ 2} \alpha'_{ 2} )  \\
\hline
1 & (1,3,2)  &  ((1,  (\alpha'_{2})^{-1}  \alpha_{ 2} \alpha'_{2} , (\alpha'_{3})^{-1}  \alpha_{ 3} \alpha'_{  3} , (\alpha'_{1})^{-1}  \alpha_{ 1} \alpha'_{ 1} )   \\
\hline
\hline
\end{array} \]
In this case the cycle decomposition of $\sigma$ is $(1)(2)(3)$ as the $j_{1}$'s run through $1,2,3$ the $k_{1}$'s run through the same set and for any triple of $\alpha'_{i}$'s we get $6$ terms each equal to 
\[ \begin{array}{l}
x^{3} \prod_{i=1}^{3}  e^{  2 \pi \sqrt{-1} {\rm Trace}_{ {\mathbb F}_{q}/ {\mathbb F}_{p}} ({\rm Trace}_{GL_{m}}(\alpha_{i} )/p    }  \\
\\
=  x^{3} e^{  2 \pi \sqrt{-1} {\rm Trace}_{ {\mathbb F}_{q}/ {\mathbb F}_{p}} ({\rm Trace}_{GL_{m}}(\alpha_{1} + \alpha_{2} )/p    } e^{  2 \pi \sqrt{-1} md/p    } . 
\end{array} \]

When $\sigma = (1,2)$ we have
\[ \begin{array}{|c|c|c|}
\hline
\hline
\sigma & \sigma' & Y^{-1}XY \\
\hline
\hline
(1,2) & 1 & ((1,2),  (\alpha'_{1})^{-1}  \alpha_{ 1} \alpha'_{ 2} ,( \alpha'_{2})^{-1}  \alpha_{ 2} \alpha'_{ 1} ,(\alpha'_{3})^{-1}  \alpha_{ 3} \alpha'_{ 3} )  \\ 
\hline
(1,2) & (1,2) &  ((1,2),  (\alpha'_{2})^{-1}  \alpha_{2} \alpha'_{ 1} , (\alpha'_{1})^{-1}  \alpha_{ 1} \alpha'_{  2} ,(\alpha'_{3})^{-1}  \alpha_{ 3)} \alpha'_{ 3} )  \\ 
\hline
(1,2) &  (1,3)  & ((2,3) ,  (\alpha'_{3})^{-1}  \alpha_{ 3} \alpha'_{  2} , (\alpha'_{2})^{-1}  \alpha_{ 2} \alpha'_{  3} ,(\alpha'_{1})^{-1}  \alpha_{ 1} \alpha'_{ 1} )  \\ 
\hline
(1,2) &  (2,3)  &  ((1,3),  (\alpha'_{1})^{-1}  \alpha_{ 1} \alpha'_{  3} , (\alpha'_{3})^{-1}  \alpha_{ 3} \alpha'_{  1} ,(\alpha'_{2})^{-1}  \alpha_{ 2} \alpha'_{  2} )  \\ 
\hline
(1,2) & (1,2,3) & ((1,3),  (\alpha'_{3})^{-1}  \alpha_{ 3} \alpha'_{ 1} , (\alpha'_{1})^{-1}  \alpha_{1} \alpha'_{  3} ,(\alpha'_{2})^{-1}  \alpha_{ 2} \alpha'_{ 2} )  \\ 
(1,2) & = (1,3,2) &  ( (2,3) ,  (\alpha'_{2})^{-1}  \alpha_{ 2} \alpha'_{ 3} , (\alpha'_{3})^{-1}  \alpha_{ 3} \alpha'_{  2} ,(\alpha'_{1})^{-1}  \alpha_{1} \alpha'_{ 1} )  \\ 
\hline
\hline
\end{array} \]

Now, if $(j_{1}, j_{2}) = (1,2) $ then we have 
\[ \begin{array}{|c|c|c|}
\hline
\hline
\sigma & \sigma' & (k_{1}, k_{2}) \\
\hline 
\hline 
(1,2) &  1 &    (1,2)    \\
\hline
(1,2) &  (1,2) &    ( 2,1)    \\
\hline
(1,2) &  (1,3)  &  (3, 2)  \\
\hline
(1,2) &  (2,3)  &  ( 1,2)  \\
\hline
(1,2) &  (1,2,3) & ( 3, 1)   \\
\hline
(1,2) &  (1,3,2) & ( 2, 3)    \\
\hline 
\hline 
\end{array} \]

From the data in these tables and the formula of Theorem \ref{7.2} we find that
\[  \begin{array}{l}
 {\rm dim}_{{\mathbb C}}(\lambda)  W_{3, G}^{x, H}(\lambda ) \\
 \\
=  \sum_{(\alpha_{1}, \alpha_{2}, 1)}  \chi_{\lambda}(1, \alpha_{1}, \alpha_{2}, 1) x^{3} 
  e^{  2 \pi \sqrt{-1} {\rm Trace}_{ {\mathbb F}_{q}/ {\mathbb F}_{p}} ({\rm Trace}_{GL_{m}}(\alpha_{1} + \alpha_{2}) )/p    }  e^{  2 \pi \sqrt{-1} md/p    }   \\
  \\
  +  \sum_{(\alpha_{1}, \alpha_{2}, 1)}   \chi_{\lambda}((1,2), \alpha_{1}, \alpha_{2}, 1)
  x^{2}  e^{  2 \pi \sqrt{-1} {\rm Trace}_{ {\mathbb F}_{q}/ {\mathbb F}_{p}} ({\rm Trace}_{GL_{m}}(\alpha_{1} \alpha_{2}) )/p    }   e^{  2 \pi \sqrt{-1} md/p    }  
  \end{array} \]
  while
  \[  \begin{array}{l}
     {\rm dim}_{{\mathbb C}}(\lambda)  W_{3, J}^{x, H}({\rm Ind}_{G}^{J}(\lambda) ) \\
  \\
  =    \sum_{  (\alpha_{1}, \alpha_{2}, 1)   } 
   \chi_{\lambda}(1 , \alpha_{1}, \alpha_{2}, 1)  x^{3}  e^{  2 \pi \sqrt{-1} {\rm Trace}_{ {\mathbb F}_{q}/ {\mathbb F}_{p}} ({\rm Trace}_{GL_{m}}(\alpha_{1} + \alpha_{2}) )/p    }  e^{  2 \pi \sqrt{-1} md/p    }  \\
   \\
 +   \frac{1}{2} \sum_{(\alpha_{1}, \alpha_{2}, 1)}   \chi_{\lambda}((1,2), \alpha_{1}, \alpha_{2}, 1)
  x^{2}  e^{  2 \pi \sqrt{-1} {\rm Trace}_{ {\mathbb F}_{q}/ {\mathbb F}_{p}} ({\rm Trace}_{GL_{m}}(\alpha_{1} \alpha_{2}) )/p    }   e^{  2 \pi \sqrt{-1} md/p    }   \\
  \\
+  \frac{1}{2}   \sum_{(\alpha_{1}, \alpha_{2}, 1)}  \chi_{\lambda}((1.2), \alpha_{1}, \alpha_{2}, 1) x^{2} 
  e^{  2 \pi \sqrt{-1} {\rm Trace}_{ {\mathbb F}_{q}/ {\mathbb F}_{p}} ({\rm Trace}_{GL_{m}}(\alpha_{1} + \alpha_{2}) )/p    }  e^{  2 \pi \sqrt{-1} md/p    }   .
\end{array} \]
\end{em}
\end{example}
\begin{remark}
\label{7.5}
\begin{em}

In Example \ref{7.4}, notice that ${\rm Ind}_{G}^{J}(\lambda)  = m( \lambda \otimes 1)$, the wreath product PSH multiplication, so that $ W_{n, \sigma_{n} \int H}^{x, H}(- )$ does not respect the PSH product.
\end{em}
\end{remark}

\section{ The fibred PSH algebra of $ GL_{n}{\mathbb F}_{q}$ and wreath products}

This section contains a couple of obvious remarks about fibred Hopf algebras\footnote{ It reflects the unhealthy
fixation with central characters, which pervades my essays \cite{Sn18} and \cite{Sn20}}.

The direct sum $R = \oplus_{n=0}^{\infty} \  R(GL_{n}{\mathbb F}_{q})$ is a PSH algebra ([23] and [22] p. 219). Here $R(GL_{n}{\mathbb F}_{q})$ is the ring of complex representations of $GL_{n}{\mathbb F}_{q}$ whose positive basis as a torsion-free abelian group consists of the irreducible representations.
Let $\hat{A} = {\rm Hom}(A, {\mathbb C}^{*})$ be the group of complex-valued characters of an abelian group $A$. Since $A$ has a multiplication map and a diagonal map the integral group ring 
${\mathbb Z}[ \hat{A}]$ is a Hopf algebra. The multiplication
\[  m : {\mathbb Z}[ \hat{A}] \otimes  {\mathbb Z}[ \hat{A}] \longrightarrow  {\mathbb Z}[ \hat{A}] \]
is the biadditive extension of the map ($\phi \in  \hat{A}, a \in A$)
\[   \phi_{1} \otimes \phi_{2} \mapsto  (a \mapsto \phi_{1}(a) \phi_{2}(a)) . \]
The comultiplication 
\[ m^{*} : {\mathbb Z}[ \hat{A}] \longrightarrow   {\mathbb Z}[ \hat{A}] \otimes  {\mathbb Z}[ \hat{A}]  \]
is given by the ${\mathbb Z}$-linear extension of $\phi \mapsto   \sum_{?? \hat{A}} \nu^{-1} \otimes   \phi \nu$.
The map which assigns to an irreducible representation its central character
gives a map of Hopf algebras
\[   \epsilon  : R \longrightarrow   {\mathbb Z}[ \hat{{\mathbb F}}_{q}^{*}]   .  \]
In fact, this map is fibred. We may write $R_{n,\phi}$ for the subgroup of $R(GL_{n}{\mathbb F}_{q}) $ generated by the irreducible representations with central character 
$\phi \in  \hat{{\mathbb F}}_{q}^{*}$.

When n = 0 we have an exceptional case where $R_{0} = {\mathbb Z}$ and there is no grading by $\phi$'s.
The multiplication 
\[ m : R(GL_{n}{\mathbb F}_{q}) \otimes R(GL_{m}{\mathbb F}_{q}) \longrightarrow  R(GL_{n+m}{\mathbb F}_{q})   \]
sends the tensor product of two positive basis elements $V \otimes W$ to
\[  {\rm Ind}_{P_{n,m}}^{GL_{n+m}{\mathbb F}_{q}} ( {\rm Inf}_{ GL_{n}{\mathbb F}_{q} \times GL_{m}{\mathbb F}_{q}}^{ P_{n,m} }(V \otimes W))  \]
where $P_{n,m}$ is the usual parabolic subgroup. Clearly, if $\phi$ and $\psi$ are the central characters of $V$ and $W$ respectively then the centre of $GL_{n+m}{\mathbb F}_{q}$ acts on this induced representation via $\phi \psi$. Hence we obtain
\[  m : R_{n,\phi} \otimes  R_{m,\psi}
 \longrightarrow  R_{n+m,\phi \psi}.   \]
The comultiplication has (a, n - a)-component given by the $U_{a,n?a}$-invariant subspace of a
$GL_{n+m}{\mathbb F}_{q}$-representation $V$ ($U_{a,n-a}$ is the unitriangular subgroup of $P_{a,n-a}$). Clearly this induces
\[    m^{*}  R_{n,\phi}  \longrightarrow  \oplus_{\nu \in \hat{{\mathbb F}}_{q}^{*} }
 \  R_{a,\nu^{-1} } \otimes  R_{n,\phi \nu} .   \]
 
 Similar remarks apply to the wreath-product PSH algebra of \S7 when, for example, $H = GL_{m}{\mathbb F}_{q}$.

 \section{fibred Shintani base change}
 
 This section contains a couple of obvious remarks about fibred Shintani base change\footnote{It reflects the unhealthy fixation Shintani base change, which pervades my essays \cite{Sn18} and \cite{Sn20}}.
 
 The Galois group ${\rm Gal}({\mathbb F}_{q^{m}}/{\mathbb F}_{q})$ acts on the PSH algebra for $GL_{n}{\mathbb F}_{q^{m}}$, which we shall denote by 
 $R({\mathbb F}_{q^{m}}) = \oplus_{n, \phi} \ R(GL_{n}{\mathbb F}_{q^{m}})_{\phi} $. If $V$ is an irreducible representation of $GL_{n}{\mathbb F}_{q^{m}}$ which is equivalent to its image under the action of ${\rm Gal}({\mathbb F}_{q^{m}}/{\mathbb F}_{q})$ then the same is true for its central character, $\phi_{V} : {\mathbb F}_{q^{m}}^{*} \longrightarrow {\mathbb C}^{*}$. By Hilbert's Theorem 90 there is a factorisation through the norm map of the form
 \[  \phi_{V} =   \lambda_{V} \cdot {\rm Norm} :  {\mathbb F}_{q^{m}}^{*} 
 \stackrel{{\rm Norm} }{\longrightarrow}  {\mathbb F}_{q}^{*} \stackrel{ \lambda_{V}}{\longrightarrow }
  {\mathbb C}^{*} .  \]
  
  The Shintani correspondence (\cite{Sn18} Chapter Nine, \S6; \cite{Sh76}) is a bijection between irreducible representations of $GL_{n}{\mathbb F}_{q}$ and Galois invariant irreducible representations of $GL_{n}{\mathbb F}_{q^{m}}$. If the central character of the Galois invariant irreducible representation $V$ is $\phi_{V}$ then the central character of the Shintani correspondent of $V$ will be $\lambda_{V}$. It is shown that the Galois invariant irreducibles generate a subalgebra of the PSH algebra and that there is an algebra isomorphism, induced by the Shintani correspondence,  between this subalgebra and the PSH algebra for ${\mathbb F}_{q}$. It is not a Hopf algebra isomorphism. 
  
  The graded subalgebra introduced above is fibred over the algebra $ {\mathbb Z}[\hat{{\mathbb F}}_{q^{m}}^{*} ] $ and the Shintani correspondence is fibred over the algebra map
 \[ {\mathbb Z}[\hat{{\mathbb F}}_{q^{m}}^{*} ]  \longrightarrow {\mathbb Z}[\hat{{\mathbb F}}_{q}^{*} ]  \]
 induced by the norm.
 
 The Gauss sums of the Shintani-related central characters are related by the Hasse-Davenport Theorem \cite{HDHH35}, which I shall recall now - for completeness. The additive measure is the character of the additive group of ${\mathbb F}_{q}$ 
 \[   \psi_{q} : {\mathbb F}_{q} \longrightarrow  {\mathbb C}^{*}  \]
 given by the composition of the trace to the prime subfield, ${\mathbb F}_{p}$ followed by
 \[  \psi_{p} :  {\mathbb F}_{p} \longrightarrow  {\mathbb C}^{*}  \]
 characterised by $\psi_{p}(z) = e^{ 2 \pi \sqrt{-1} z/p}$.
 
 If $\lambda : {\mathbb F}_{q}^{*} \longrightarrow  {\mathbb C}^{*} $ is a character then its Gause sum is defined by
 \[    \tau( \lambda ) = \sum_{x \in {\mathbb F}_{q}^{*} } \ \phi(x) \psi_{q}(x) .\]
 If $\phi : {\mathbb F}_{q^{m}}^{*} \longrightarrow  {\mathbb C}^{*} $ is a character which factorises through the norm as
 \[  {\mathbb F}_{q^{m}}^{*}  \stackrel{{\rm Norm}}{\longrightarrow}  {\mathbb F}_{q}^{*} \stackrel{\lambda}{\longrightarrow}   {\mathbb C}^{*}  \]
 then
 \[   -  \tau( \phi) =   (-1)^{m} \tau(\lambda)^{m} . \]

  \section{The hyperHecke algebra of a locally profinite group }
  
  This section is a recapitulation from \cite{Sn20}.
 
 Let $G$ be a locally profinite group and let $k$ be an algebraically closed field. Suppose that $\underline{\phi} : Z(G) \longrightarrow k^{*}$ is a fixed $k$-valued, continuous character on the centre $Z(G)$ of $G$. Let ${\mathcal M}_{cmc, \underline{\phi}}(G)$ be the poset of pairs $(H, \phi)$ where $H$ is a subgroup of $G$, containing $Z(G)$, which is compact, open modulo the centre of $G$ and $\phi : H \longrightarrow k^{*}$ is a $k$-valued, continuous character whose restriction to $Z(G)$ is $\underline{\phi}$.
 
 We define the hyperHecke algebra, ${\mathcal H}_{cmc}(G)$,  to be the $k$-algebra given by the following generators and relations. For $(H, \phi), (K, \psi) \in {\mathcal M}_{cmc,\underline{\phi}}(G)$,  write 
$[(K, \psi), g, (H, \phi)]$ for any triple consisting of $g \in G$, characters $\phi, \psi$ on 
subgroups $H, K \leq G$, respectively such that 
\[    (K, \psi) \leq (g^{-1}Hg, (g)^{*}(\phi)) \]
which means that $K \leq  g^{-1}Hg$ and that $\psi(k) = \phi(h)$ where $k = g^{-1}hg$ for
 $h \in H, k \in K$.

 Let ${\mathcal H}$ denote the $k$-vector space with basis given by these triples. Define a product on these triples by the formula 
\[  [(H, \phi), g_{1}, (J, \mu)]  \cdot  [(K, \psi), g_{2}, (H, \phi)] =   [(K, \psi), g_{1}g_{2}, (J, \mu)]  \]
and zero otherwise. This product makes sense because 

(i)  \   if $K \leq g_{2}^{-1} H g_{2}$ and
 $H \leq g_{1}^{-1} J g_{1}$ then $K \leq  g_{2}^{-1} H g_{2} \leq   g_{2}^{-1} g_{1}^{-1} J g_{1} g_{2} $ 
 
 and 
 
 (ii)  \  if $\psi(k) = \phi(h) = \mu(j), $ where $k = g_{2}^{-1}hg_{2},  h = g_{1}^{-1}j g_{1}$ then
 \linebreak
  $k = g_{2}^{-1}  g_{1}^{-1}j g_{1}   g_{2}$. 
 
 This product is clearly associative and we define an algebra ${\mathcal H}_{cmc}(G)$ to be ${\mathcal H}$ modulo the relations ( \cite{Sn18}\footnote{For the purposes of this essay, as in \cite{Sn20},  I  am using throughout a convention in which $g$ is replaced by $g^{-1}$ to make the composition coincide with the conventions of induced representations which I shall use here.})
 \[   [(K, \psi), gk, (H, \phi)]  = \psi(k^{-1}) [(K, \psi), g, (H, \phi)]  \]
and
 \[     [(K, \psi), hg, (H, \phi)]  = \phi(h^{-1}) [(K, \psi), g, (H, \phi)] .    \]

\section{A Hopf-like algebra for the hyperHecke algebras of $GL_{n}K$}

 Let $G$ be a locally profinite group and let $k$ be an algebraically closed field. A character of a group $H$ will mean a continuous homomorphism $\phi : H \longrightarrow k^{*}$. Usually it is important that pairs $(H, \pi)$ belong to the poset ${\mathcal M}_{cmc,\underline{\phi}}(G)$, introduced in \S10. 
 
 However, the material of this section, while applying without change to the case of general linear groups of local fields, is amply illustrated by the case of general linear groups of a finite field. In particular, the finite field case allows us to use the tensor product description of  induced representations rather than
the function space version\footnote{These are compared in one of the Appendices to \cite{Sn20}.}. 

Recall that the hyperHecke algebra of $G$ is denoted by ${\mathcal H}_{cmc}(G)$\footnote{The subscript stands for ``compact open modulo the centre''.}.

As explained in \S10  ${\mathcal H}_{cmc}(G)$ is generated by triples 
$[(K, \psi) , g, (H, \phi)]$ where $H,K \subseteq G$ and $\phi, \psi$ are characters such that $(K, \psi) \leq  (g^{-1}Hg, g^{*}(\phi))$. This condition means that $K \subseteq g^{-1}Hg$ and if $k \in K$ satisfies $gkg^{-1} = h \in H$ then $\psi(k) = \phi(h)$.

We think, as motivation, of a triple as a $k[G]$-module homomorphism
\[  [(K, \psi) , g, (H, \phi)]  : {\rm Ind}_{K}^{G}(\psi) \longrightarrow {\rm Ind}_{H}^{G}(\phi) \]
defined, in the tensor product description of induced representations (see \cite{Sn20} Appendix \S12), by $g' \otimes_{K} v \mapsto g' g^{-1} \otimes_{H} v$. 

This motivation with inform the constructions which follow. For 
\linebreak
$G_{1} = GL_{n}{\mathbb F}_{q}$ and $G_{2} = GL_{m}{\mathbb F}_{q}$
suppose we have  $[(K_{1}, \psi_{1}) , g_{1}, (H_{1}, \phi_{1})] $ and 
\linebreak
$[(K_{2}, \psi_{2}) , g_{2}, (H_{2}, \phi_{2})] $  in ${\mathcal H}_{cmc}(G_{1})$ and ${\mathcal H}_{cmc}(G_{2})$, respectively.

Inside $GL_{n+m}{\mathbb F}_{q}$ we have the unitriangular subgroup $U_{n,m}$ such that the parabolic subgroup $P_{n,m} = (G_{1} \times G_{2})U_{n,m}$. On $P_{n,m}$ we have ${\rm Inf}_{K_{1} \times K_{2}}^{K_{n,m}}(\psi_{1} \otimes \psi_{2})$ and 
${\rm Inf}_{H_{1} \times H_{2}}^{H_{n,m}}(\phi_{1} \otimes \phi_{2})$ where
$K_{n,m} = (K_{1} \times K_{2})U_{n,m}$ and 
\linebreak
$H_{n,m} = (H_{1} \times H_{2})U_{n,m}$.

Define the product of these two triples, as a triple in ${\mathcal H}_{cmc}(GL_{n+m}{\mathbb F}_{q})$, to equal
\[ \begin{array}{l}
m( [(K_{1}, \psi_{1}) , g_{1}, (H_{1}, \phi_{1})]  \otimes  [(K_{2}, \psi_{2}) , g_{2}, (H_{2}, \phi_{2})] \\
\\
= [ (K_{n,m}, {\rm Inf}(\psi_{1} \otimes \psi_{2})), (g_{1}, g_{2})  , (H_{n,m}, {\rm Inf}(\phi_{1} \otimes \phi_{2})]   .
  \end{array} \]
  
  Next we need to motivate the coproduct by recalling the double coset formula (cite{Sn94} Theorem 1.2.40; \cite{Sn18} Chapter Seven).

If $K \subseteq GL_{n+m}{\mathbb F}_{q}$ we have 
  \[ \begin{array}{l} 
  {\rm Res}_{P_{a, n+m-a}}^{GL_{n+m}{\mathbb F}_{q}} ({\rm Ind}_{K}^{GL_{n+m}{\mathbb F}_{q}}( \psi)  ) \\
  \\
\cong  \oplus_{z \in P_{a,n+m-a} \backslash GL_{n+m}{\mathbb F}_{q} / K }  \
  {\rm Ind}_{P_{a,n+m-a} \bigcap  z K z^{-1} }^{P_{a,n+m-a}}( (z^{-1})^{*}( \psi)).
  \end{array} \]
  The explicit maps which give the Double Coset Formula in this case are
  \[  g \otimes_{K} w   \mapsto  u \otimes_{P_{a,n+m-a} \bigcap  z K z^{-1}} kw , \ g = uzk, u \in P_{a,n+m-a} , k \in K       \]
  and 
\[ u \otimes_{P_{a,n+m-a} \bigcap  z K z^{-1}} v \mapsto  uz \otimes_{K} v .\] 

Notice that 
\[ g \otimes_{K} w   \mapsto  u \otimes_{P_{a,n+m-a} \bigcap  z K z^{-1}} kw \mapsto
uz \otimes_{ K} kw = uzk \otimes_{ K} w =   g \otimes_{K} w .\]
Also
\[ u \otimes_{P_{a,n+m-a} \bigcap  z K z^{-1}} w \mapsto  uz \otimes_{ K} w  \mapsto 
u \otimes_{P_{a,n+m-a} \bigcap  z K z^{-1}} w  \]
so if the maps are well-defined then they are inverse isomorphisms.

I shall pause to check that they are well-defined, as this process may be useful in what follows.

If $g = uzk = u' z' k'$ then  $gKg^{-1} = uzKz^{-1}u^{-1} = u' z'K(z')^{-1} (u')^{-1}$ so we have
\[   (u^{-1}u' -  (u')^{-1}u ) :     P_{a,n+m-a} \bigcap  z'K(z')^{-1} \stackrel{\cong}{\longrightarrow}   P_{a,n+m-a}  \bigcap zKz^{-1}  \]
also if $v \in P_{a,n+m-a} \bigcap  z'K(z')^{-1}$ it acts via $\psi((z')^{-1}vz')$ whereas 
$(u^{-1}u' v  (u')^{-1}u )$ acts via 
\[ \begin{array}{l}
 \psi(z^{-1}u^{-1}u' v  (u')^{-1}uz) = \psi(k^{-1}z^{-1}u^{-1}u' v  (u')^{-1}uzk)   \\
 \\
 = \psi( k^{-1} (z')^{-1}   v z'k) = \psi( (z')^{-1}   v z') .
 \end{array}  \]

 Recall that $(K, \psi) \leq  (g^{-1}Hg, g^{*}(\phi))$ means that\footnote{Apologies for $k$ the field and $k \in K$. The former will hardly appear again!} 
 \[     K \stackrel{\psi}{\rightarrow} {\mathbb C}^{*} =  K \stackrel{(g - g^{-1})}{\rightarrow} H 
 \stackrel{\phi}{\rightarrow}  k^{*} \]
 or equivalently $k = g^{-1}hg$ implies $\psi(k) = \phi(h)$.
 
 If $\underline{K} =  ( P_{a,n+m-a} \bigcap  z' K (z')^{-1}$ and $\underline{H} = P_{a,n+m-a} \bigcap  z K z^{-1}$ and $\underline{g} = u^{-1}u'$ we have
  \[     \underline{K} \stackrel{\psi( (z')^{-1}   - z')}{\rightarrow} {\mathbb C}^{*} =  \underline{K} \stackrel{(\underline{g} - \underline{g}^{-1})}{\rightarrow} \underline{H} 
 \stackrel{\psi(z^{-1} - z)}{\rightarrow}  k^{*} \]
 
 From the first datum we get
 \[  [(K, \psi) , g, (H, \phi)]  : {\rm Ind}_{K}^{G}(\psi) \longrightarrow {\rm Ind}_{H}^{G}(\phi) \]
defined by $g' \otimes_{K} v \mapsto g' g^{-1} \otimes_{H} v$.

Therefore from the second we get ($P = P_{a,n+m-a}$)
 \[  \begin{array}{l}
  [(\underline{K}, \psi((z')^{-1} - z')) , \underline{g}, (\underline{H}, \psi( (z)^{-1}   - z))]  : 
 {\rm Ind}_{\underline{K}}^{P}(\psi((z')^{-1} - z')) \\
 \\
 \hspace{150pt} 
 \stackrel{\cong}{ \longrightarrow} {\rm Ind}_{\underline{H}}^{P}(\psi(z^{-1} - z)) 
 \end{array} \]
defined by 
\[ g' \otimes_{\underline{K}} v \mapsto g' \underline{g}^{-1} \otimes_{\underline{H}} v 
=  g' (u^{-1}u')^{-1} \otimes_{\underline{H}} v =   g' (u')^{-1}u \otimes_{\underline{H}} v   .  
  \]
 
 The double-check above confirms that conjugation as above on $P_{a,n+m-a}$ gives an canonical isomorphism 
which, in the language of hyperHecke triples gives an inclusion (an isomorphism) 
 \[      \begin{array}{l}
        ( P_{a,n+m-a} \bigcap  z' K (z')^{-1} , ((z')^{-1})^{*}( \psi))    \\
        \\
        \hspace{80pt}                         \leq   (  (u^{-1}u' )^{-1}P_{a,n+m-a} \bigcap  z K z^{-1}u^{-1}u'  , ((u^{-1}u' )^{*} (z^{-1})^{*}( \psi))  .
        \end{array}  \]
 Therefore, to recapitulate,  we have a triple
 \[ [   ( P_{a,n+m-a} \bigcap  z' K (z')^{-1} , ((z')^{-1})^{*}( \psi)) , u^{-1}u' ,  
P_{a,n+m-a} \bigcap  z K z^{-1}  , (z^{-1})^{*}( \psi))  ] \]
which gives a canonical $P_{a,n+m-a}$-isomorphism of induced modules 
\[  \begin{array}{l}
{\rm Ind}_{P_{a,n+m-a} \bigcap  z' K (z')^{-1} }^{P_{a,n+m-a}}(  ((z')^{-1})^{*}( \psi)) ) \longrightarrow 
{\rm Ind}_{P_{a,n+m-a} \bigcap  z K (z)^{-1} }^{P_{a,n+m-a}}(  ((z)^{-1})^{*}( \psi)) )
\end{array} \]
given by 
\[    u_{1} \otimes_{P_{a,n+m-a} \bigcap  z' K (z')^{-1} } v  \mapsto 
u_{1}(u^{-1}u' )^{-1} \otimes_{P_{a,n+m-a} \bigcap  z K (z)^{-1} }    v .   \]

Now, returning to the construction of the coproduct, suppose we have 
\[ [(K, \psi) , g , (H, \phi)] : {\rm Ind}_{K}^{GL_{n+m}{\mathbb F}_{q}}( \psi)  \longrightarrow  {\rm Ind}_{H}^{GL_{n+m}{\mathbb F}_{q}}( \phi)  \]
given by $ g' \otimes_{K} v \mapsto g' g^{-1} \otimes_{H} v$ with $g' = uzk$, $u \in P_{a.n+m-a}, k \in K, z \in G$. This is a $GL_{n+m}{\mathbb F}_{q}$ map so it is a $P_{a,n+m-a}$-map and therefore will
induce a lot of maps between the summands in the Double Coset Formula.

The left-hand map to the Double Coset summand is, if $g' = uzk$ as above, 
\[    g' \otimes_{K} v  \longrightarrow  u \otimes_{P_{a,n+m-a} \bigcap zKz^{-1} } kv  \]
and the right-hand map to  the Double Coset summand requires us to find the double coset of 
$g' g^{-1} = u z k g^{-1} = u zg^{-1} (gkg^{-1}) \in P_{a,n+m-a} zg^{-1} H$ therefore the right-hand map is
\[   g' g^{-1} \otimes_{H} v  \longrightarrow  u \otimes_{P_{a,n+m-a} \bigcap zg^{-1}Hgz^{-1} } (gkg^{-1})v \]
Therefore the map between Double Coset summands corresponding to 
\linebreak
$[(K, \psi) , g , (H, \phi)]$
is
\[   u \otimes_{P_{a,n+m-a} \bigcap zKz^{-1} } kv \longrightarrow  u \otimes_{P_{a,n+m-a} \bigcap zg^{-1}Hgz^{-1} } (gkg^{-1})v \]
but if we insert the $K$ and $H$ actions in terms of $\psi$ and $\phi$ it becomes
\[  \begin{array}{l}
u \otimes_{P_{a,n+m-a} \bigcap zKz^{-1} } \psi(k) v \longrightarrow  u \otimes_{P_{a,n+m-a} \bigcap zg^{-1}Hgz^{-1} } \phi(gkg^{-1})v \\
\\
\hspace{180pt} =  u \otimes_{P_{a,n+m-a} \bigcap zg^{-1}Hgz^{-1} } \psi(k)v 
\end{array}  \]
or simply
\[  u \otimes_{P_{a,n+m-a} \bigcap zKz^{-1} }  v \longrightarrow    u \otimes_{P_{a,n+m-a} \bigcap zg^{-1}Hgz^{-1} } v . \]
Since
\[ (z^{-1})^{*}(\psi) :  P_{a,n+m-a} \bigcap zKz^{-1} \longrightarrow k^{*} \]
equals
\[ P_{a,n+m-a} \bigcap zKz^{-1} \leq P_{a,n+m-a} \bigcap zg^{-1}Hgz^{-1}
\stackrel{ ((zg^{-1})^{-1})^{*}(\phi)}{\longrightarrow} k^{*} \]
the map between Double Coset summands is
\[ [ ( P_{a,n+m-a} \bigcap zKz^{-1} ,  (z^{-1})^{*}(\psi) ) ,  1  , (P_{a,n+m-a} \bigcap zg^{-1}Hgz^{-1} , 
((zg^{-1})^{-1})^{*}(\phi) ) ]   \]
between induced representations of the form
\[ \begin{array}{c}
  {\rm Ind}_{P_{a,n+m-a} \bigcap zKz^{-1}}^{P_{a,n+m-a}}(  (z^{-1})^{*}(\psi) )  \\
  \\
  \downarrow    \\
  \\
  {\rm Ind}_{P_{a,n+m-a} \bigcap zg^{-1}Hgz^{-1}}^{P_{a,n+m-a}}(  ((zg^{-1})^{-1})^{*}(\phi)    ) .
\end{array} \]
For each triple we have one such $P_{a,n+m-a}$-map for each $0 \leq a \leq n+m$.

Finally, from this triple, we want to derive a triple lying in ${\mathcal H}_{cmc}(GL_{a}{\mathbb F}_{q} \times GL_{n+m-a}{\mathbb F}_{q} )$. This will use the isomorphism
\[ P_{a,n+m-a}/U_{a,n+m-a} \cong  GL_{a}{\mathbb F}_{q}  \times GL_{n+m-a}{\mathbb F}_{q} . \]
In the world of representations one would just average over the $U_{a,n+m-a}$-action but we are in a much more combinatorial situation. 

Instead, suppose that the character $((zg^{-1})^{-1})^{*}(\phi) =1$ on 
\linebreak
$U_{a,n+m-a} \bigcap zg^{-1}Hgz^{-1}$, which guarantees that $((z^{-1})^{-1})^{*}(\psi) =1$ on 
\linebreak
$U_{a,n+m-a} \bigcap z^{-1}Kz^{-1}$. This is true for all when $a=0$ or $a=n+m$.

 In this case our hyperHecke triple induces a triple in the hyperHecke algebra of $(P_{a,n+m-a}/U_{a,n+m-a}) \cong  GL_{a}{\mathbb F}_{q} \times GL_{n+m-a}{\mathbb F}_{q}$
\[ [ ( \frac{P_{a,n+m-a} \bigcap zKz^{-1}}{U_{a,n+m-a} \bigcap zKz^{-1}} ,  (z^{-1})^{*}(\psi) ) ,  1  ,
 ( \frac{ P_{a,n+m-a} \bigcap zg^{-1}Hgz^{-1}}{U_{a,n+m-a} \bigcap zg^{-1}Hgz^{-1} } , 
((zg^{-1})^{-1})^{*}(\phi) ) ]   \]
which lies in 
\[ {\mathcal H}_{cmc}(GL_{a}{\mathbb F}_{q} \times GL_{n+m-a}{\mathbb F}_{q} ) \cong 
 {\mathcal H}_{cmc}(GL_{a}{\mathbb F}_{q}  ) \otimes  {\mathcal H}_{cmc}( GL_{n+m-a}{\mathbb F}_{q} ) .\]

This construction works equally well for $GL_{n}$ of a local field.
\begin{definition}{(hyperHecke ``coproduct'')}
\label{11.1}
\begin{em}
 Let $F$ be a non-Archimedean local field or a finite field. The ``coproduct''
 \[ m^{*} :  {\mathcal H}_{cmc}(GL_{n}F)  \longrightarrow  \oplus_{a=0}^{n} \ {\mathcal H}_{cmc}(GL_{a}F \times  GL_{n-a}F)  \]
 is defined by the formula
 \[ \begin{array}{l}
m^{*}([(K, \psi) , g , (H, \phi)] )_{a,n+m-a} \\
\\
= \tilde{\sum}_{a=0}^{n}  \  [ ( \frac{P_{a,n+m-a} \bigcap zKz^{-1}}{U_{a,n+m-a} \bigcap zKz^{-1}} ,  (z^{-1})^{*}(\psi) ) ,  1  ,
 ( \frac{ P_{a,n-a} \bigcap zg^{-1}Hgz^{-1}}{U_{a,n+m-a} \bigcap zg^{-1}Hgz^{-1} } , 
((zg^{-1})^{-1})^{*}(\phi) ) ]  
\end{array}  \]
where $\tilde{\sum}_{a=0}^{n} $ indicates that only the terms for which 
$((zg^{-1})^{-1})^{*}(\phi) =1$ on $U_{a,n-a} \bigcap zg^{-1}Hgz^{-1}$ are summed over.
\end{em}
\end{definition} 

\begin{remark}{Hopflike algebra}
\label{11.2}
\begin{em}

 The purpose of this remark is to clarify the notion of a Hopflike algebra.
\end{em}
\end{remark} 

Definition \ref{11.1} constructs a ``coproduct'' 
 \[  m^{*} :  {\mathcal H}_{cmc}(GL_{n}F)  \longrightarrow  \oplus_{a=0}^{n} \ {\mathcal H}_{cmc}(GL_{a}F \times  GL_{n-a}F)  \] 
where $ {\mathcal H}_{cmc}(G) $ is the hyperHecke algebra of a locally $p$-adic Lie group $G$.

I believe the coproduct, correctly interpeted, is an ``algebra map'' with respect to the multiplication $m$. for example the calculations of \cite{Sn20c} show that
\[   {\mathcal H}_{cmc}(GL_{a}F)  \otimes  {\mathcal H}_{cmc}(GL_{n-a}F) 
\stackrel{m}{\longrightarrow}  {\mathcal H}_{cmc}(GL_{n}F) \stackrel{m^{*}}{\longrightarrow}  {\mathcal H}_{cmc}(GL_{b}F  \times  GL_{n-b}F) \]
and 
\[   \begin{array}{l}
 {\mathcal H}_{cmc}(GL_{a}F)  \otimes  {\mathcal H}_{cmc}(GL_{n-a}F) \\
 \\
\stackrel{m^{*} \otimes m^{*}}{\longrightarrow}  
 {\mathcal H}_{cmc}(GL_{x_{1,1}}F  \times GL_{x_{1,2}}F \times GL_{x_{2,1}}F  \times GL_{x_{2,2}}F) \\
 \\
 \stackrel{1 \times T \times 1}{\longrightarrow}
  {\mathcal H}_{cmc}(GL_{x_{1,1}}F  \times GL_{x_{2,1}}F \times GL_{x_{1,2}}F  \times GL_{x_{2,2}}F)
  \\
  \\
    \stackrel{m \times m}{\longrightarrow}  {\mathcal H}_{cmc}(GL_{b}F  \times  GL_{n-b}F) 
    \end{array}  \]
    are equal for any matrix of positive integers 
    \[   \left( \begin{array}{cc}
    x_{1,1} & x_{1,2} \\
    \\
    x_{2,1} & x_{2,2} 
    \end{array} \right)\]
 satisfying 
 \[ x_{1,1} + x_{1,2} = a, x_{2,1}+ x_{2,2} = n-a, x_{1,1} + x_{2,1} = b, x_{1,2} + x_{2,2} = n-b .\]
The calculation uses the Bruhat decomposition in the combinatorial manner of (\cite{AVZ81} Appendix III).

It seems to me that there is a family of analogous pairs of equal compositions indexed by ordered partitions
 $\underline{n} = (n_{1}, n_{2}, \ldots , n_{t})$ and featuring $ {\mathcal H}_{cmc}(GL_{n_{1}}F \times \ldots \times GL_{n_{t}}F)$, each of which is a consequence of the Bruhat decomposition. The one given above and  proved in \cite{Sn20c} is the case of partitions of length two.
 
 I conjecture that the result is a Hopf algebra structure on the sum of the hyperHecke algebras of products of general linear groups, indexed by ordered partitions.

\end{document}